\documentclass[a4paper,11pt]{amsart}

\usepackage{amsmath}
\usepackage{amssymb}
\usepackage{amsthm}
\usepackage[all]{xy}
\usepackage{textcomp}

\newtheorem{lemma}{Lemma}[section]

\newtheorem{thm}{Theorem}[section]

\newcommand{\G}{\mathrm G}
\newcommand{\pG}{\tilde{\mathrm{G}}}
\newcommand{\gK}{\mathrm K}
\newcommand{\gL}{\mathrm L}
\newcommand{\gH}{\mathrm H}

\newcommand{\mV}{\mathrm V}

\newcommand{\mW}{\mathrm W}
\newcommand{\mZ}{\mathrm Z}
\newcommand{\SO}{\mathrm{SO}}
\newcommand{\CSO}{\mathrm{CSO}}
\newcommand{\CO}{\mathrm{CO}}
\newcommand{\gO}{\mathrm O}

\newcommand{\gU}{\mathrm U}
\newcommand{\SL}{\mathrm{SL}}
\newcommand{\End}{\mathrm{End}}
\newcommand{\Sp}{\mathrm{Sp}}
\newcommand{\pSp}{\mathrm{\widetilde{Sp}}}
\newcommand{\CSp}{\mathrm{CSp}}
\newcommand{\GL}{\mathrm{GL}}
\newcommand{\Spin}{\mathrm{Spin}}

\newcommand{\gM}{\mathrm M}
\newcommand{\gR}{\mathrm{R}}
\newcommand{\tr}{\mathrm{tr}}
\newcommand{\Stab}{\mathrm{Stab}}

\newcommand{\C}{\mathbb C}
\newcommand{\R}{\mathbb R}
\newcommand{\Q}{\mathbb H}
\newcommand{\Oc}{\mathbb O}
\newcommand{\pOc}{{\tilde{\mathbb O}}}
\newcommand{\pQ}{{\tilde{\mathbb H}}}
\newcommand{\Z}{\mathbb Z}

\newcommand{\cB}{\mathcal{B}}

\newcommand{\p}{\tilde}
\newcommand{\wt}{\widetilde}

\newcommand{\CD}{\mathfrak{CD}}
\newcommand{\pCD}{\widetilde{\CD}}

\newcommand{\pI}{{\tilde I}}
\newcommand{\pJ}{{\tilde J}}
\newcommand{\pK}{{\tilde K}}

\newcommand{\ra}{\rightarrow}

\title{Multisymplectic 3-forms on 7-dimensional manifolds}
\author{Tom\'a\v s Sala\v c}
\begin{document}

\maketitle

\begin{abstract}
A 3-form $\omega\in\Lambda^3\R^{7\ast}$ is called multisymplectic if it satisfies some natural non-degeneracy requirement. It is well known that there are  8 orbits (or types) of multisymplectic 3-forms on $\R^7$ under the canonical action of $\GL(7,\R)$ and that two types are open. This leads to 8 types of global multisymplectic 3-forms on 7-dimensional manifolds without boundary. 
The existence of a global multisymplectic 3-form of a fixed type is a classical problem in differential topology which is equivalent to the existence of a certain $G$-structure. The open types are the most interesting cases as they are equivalent  to a $\G_2$ and  $\pG_2$-structure, respectively. The existence of these two structures is a well known and solved problem. In this article is solved (under some convenient assumptions) the problem of the existence multisymplectic 3-forms  of the remaining  types.  
\end{abstract}

\section{Introduction}
Put  $\mV:=\R^7$. 
There are finitely many orbits of  the canonical action of $\GL(\mV)$  on $\Lambda^3\mV^\ast$. 
We will call the orbits also types. A linear isomorphism $\Phi:\R^7\ra\mW$ induces a map $\Phi^\ast:\Lambda^3\mW^\ast\ra\Lambda^3\R^{7\ast}$. The type of $\Phi^\ast\omega$ does not depend on the choice of linear isomorphism and thus, we can define the type for any skew-symmetric 3-form on any 7-dimensional real vector space. 

A 3-form $\omega\in\Lambda^3\mV^\ast$ is called \textit{multisymplectic} if the insertion map 
\begin{equation}
\mV\ra\Lambda^2\mV^\ast, \ v\mapsto i_v\omega:=\omega(v,-,-)
\end{equation}
is injective. There are (see \cite{Dj} and \cite{W}) eight types of multisympletic 3-forms  and  two open types.

Let $\Omega$ be a global 3-form on a 7-dimensional manifold $N$ without boundary and $i\in\{1,\dots,8\}$. We call $\Omega$ a \textit{multisympletic 3-form of algebraic type}  $i$ if  for each $x\in N:\ \Omega_x$ is   a multisymplectic 3-form of type $i$.  The existence of such a 3-form is a classical problem in differential topology, if $\gO_i$ is the stabilizer of a fixed multisymplectic 3-form $\omega_i\in\Lambda^3\mV^\ast$  of algebraic type $i$, then $N$ admits a multisymplectic 3-form of algebraic type $i$ if, and only if it has an $\gO_i$-structure.  The groups $\gO_i$ were studied in \cite{BV} where they were given as semi-direct products of some well known Lie groups. 

By the Cartan-Iwasawa-Malcev theorem (see \cite[Theorem 1.2]{Bo}), a connected Lie group $\gH$ has a maximal compact subgroup and any two such subgroups are conjugated. Let us fix one such subgroup and let us denote it by $\gK$. By Cartan's result, the group $\gH$ has the homotopy type of $\gK$ and by a standard argument from the obstruction theory, any $\gH$-principal bundle  reduces to a $\gK$-principal bundle.
Hence, the first goal  is (see Section \ref{section max compact subgroups}) to find a maximal compact subgroup $\gK_i$ of each group $\gO_i$. Then  we solve (see Section \ref{section global forms}) the problem of the existence of a multisympletic form on a closed 7-manifold of algebraic type $i$. The problem is not solved completely as for some types we assume that the underlying manifold is orientable or simply-connected.

The most interesting and well known cases are types 8 and 5 as $\gO_8=\G_2$ and $\gO_5=\pG_2$. The existence of a $\G_2$-structure was solved in \cite{G} and the existence of a  $\pG_2$-structure  in \cite{Le}. 

\bigskip

Let us summarize the main result of this article into a single Proposition. See Section \ref{section spin char class} for the definition of  characteristic classes $q(N)$ and $q(N;\ell)$.

\begin{thm}\label{metatheorem}
Let $N$ be a closed and connected 7-manifold. 
\begin{enumerate}
 \item Suppose that $N$ is orientable, spin$^c$ and that there are $e,f\in H^2(N,\Z)$ such that 
 \begin{equation*}
  w_2(N)=\rho_2(e+f)\ \ \mathrm{and}\ \ q(N;e+f)=-ef,
 \end{equation*}
then $N$ admits a multisymplectic 3-form  of algebraic type 1.

If $N$ is simply-connected, then the assumptions are also necessary.
 \smallskip
\item  Suppose that $N$ is orientable, spin and that there are  $e,f\in H^2(N,\Z)$ such that  $$-q(N)=e^2+f^2+3ef,$$ 
then $N$ admits a multisymplectic 3-form  of algebraic type 2.

If $N$ is simply-connected, then the assumptions are also necessary.
\smallskip
\item $N$ admits a multisympletic 3-form of algebraic type $3$ if, and only if $N$ is orientable and spin$^c$.
\smallskip
\item  Suppose that  $N$ is orientable, spin and there is $u\in H^4(N,\Z)$ such that $$q(M)=-4u,$$
then $N$ admits a multisymplectic 3-form  of algebraic type 4.

On the other hand, if $N$ admits a multisymplectic 3-form of algebraic type $4$, then $N$ is orientable and spin.
\smallskip
\item $N$ admits a multisympletic 3-form of algebraic type $i=5,6,7,8$ if, and only if $N$ is orientable and spin.
\end{enumerate}
\end{thm}

\bigskip

\textit{Acknowledgement and dedication}.  The author is grateful to Micheal C. Crabb for pointing out several mistakes in the original article and considerably simplifying many arguments in the present work. I would also like to thank to J. Van\v zura, M. \v Cadek  and to the unknown referee for several valuable comments and suggestions. 

\bigskip

I would like to dedicate this article to M. Doubek who drew author's attention to this subject and who passed away in a car accident at the  age of 33.

\subsection{Notation}\label{section notation}
We will use  the following notation: 
\bigskip

\begin{tabular}{rl}
$1_n:=$& identity $n\times n$ matrix,\\
$[v_1,\dots,v_i]:=$& the linear  span of vectors $v_1,\dots,v_i$,\\
$\alpha_{i_1i_2\dots i_\ell}:=$&$\alpha_{i_1}\wedge\alpha_{i_2}\wedge\dots\wedge\alpha_{i_\ell}$,\\
$M(k,\R):=$& the algebra of $k\times k$ real matrices.
\end{tabular}

\bigskip

Let $\mV,\mW$ be real vector spaces, $\mV^\ast$ be the dual vector space to $\mV$, $\End(\mV)$ be the algebra of linear endomorphisms of $\mV$, $\GL(\mV)$ be the group of linear automorphisms of $\mV$. Suppose that $\mW'$ is a vector subspace of $\mW$,\ \ $A$ is a subalgebra of $\End(\mV)$,\ \ $\omega\in\otimes^i\mW^\ast$,\ \  $\varphi\in\GL(\mW)$,\ \ $H\subset\GL(\mW)$, \ \ $ \cB=\{v_1,\dots,v_n\}$ is a basis of $\mV$ and  $\Phi:\mV\ra\mW$ is a linear isomorphism. Then we put 
\begin{align}
&\Phi^\ast\varphi:=\Phi^{-1}\circ\varphi\circ\Phi\in\GL(\mV)\ \ \mathrm{ and} \ \ \varphi^\ast H:=\{\varphi^\ast h:\ h\in H\},\label{notation pullback of map}\\
&\Phi^\ast\omega\in\otimes^i\mV^\ast,\ \Phi^\ast\omega(v_1,\dots,v_i):=\omega(\Phi(v_1),\dots,\Phi(v_i)), \ v_1,\dots,v_i\in\mV,\\\label{notation pullback of tensor}
&\varphi.\omega:=(\varphi^{-1})^\ast\omega,\\
&\Stab(\omega):=\{\varphi\in\GL(\mW):\ \varphi.\omega=\omega\},\label{notation stabilizer of tensor}\\
&\Stab(\mW'):=\{\varphi\in\GL(\mW):\ \varphi(\mW')\subset\mW'\},\label{notation stabilizer of subspace}\\
&\GL_A(\mV):=\{\varphi\in\GL(\mV):\ \forall a\in A,\ \varphi\circ a=a\circ\varphi\},\label{notation equivariant maps}\\
&\phi_{\cB}:=(a_{ij})\in M(n,\R),\ \phi\in\End(\mV),\ \phi(v_j)=\sum_i a_{ij}v_i,\\
&\triangle^k(\omega):=\{w\in\mW:\ (i_w\omega)^{\wedge^k}=0\}\label{notation isotropy set}.
\end{align}
It is easy to see that $\varphi^\ast\Stab(\omega)=\Stab(\varphi^\ast\omega)$.
\bigskip

From now on, $\cB_{st}=\{e_1,\dotsm,e_\ell\}$ will denote the standard basis of $\R^\ell,\ \ell\ge1$. If $f_1,\dots,f_\ell\in\R^\ell$, then $(f_1|\dots|f_\ell)$ denotes the matrix with $i$-th column equal to $f_i$. We put $E_{ij}:=(0|\dots|e_i|\dots|0)\in M(\ell,\R)$ where $e_i$ is 
the $j$-th column. We for brevity put $\mV:=\R^7$ and  $\triangle^j_i:=\triangle^j(\omega_i)$ where  $\omega_i$ is the fixed  multisymplectic 3-form of type $i$ from \cite{BV}. It follows that $\triangle^j_i$ is an $\gO_i$-invariant subset.

\section{$\ast$-Algebras, some subgroups of $\G_2,\tilde\G_2$ and semi-direct products of Lie groups}
In this Section we will review some well known theory, set notation and prove some elementary facts about maximal compact subgroups of semi-direct products of Lie groups. 

\subsection{$\ast$-algebras}\label{section algebras}

A $\ast$-\textit{algebra} is a pair $(A,\ast)$ where:
\begin{itemize}
 \item $A$ is a  real algebra\footnote{We do not assume that $A$ is associative.} with a unit $1$ and 
 \item  $\ast:A\ra A$ is an involution  which satisfies: $$\ast 1=1,\ \ast(a.b)=(\ast b).(\ast a), \ a,b\in A.$$
\end{itemize}
We put  $\bar a:=\ast(a)$,  $\Re A:=\{a\in A:\ \bar a=a\}$ and $\Im A:=\{a\in A:\ \bar a=-a\}$. Then  $A=\Re A\oplus\Im A$ and we denote by $\Re(a)$ and $\Im(a)$ the projection of  $a\in A$  onto $\Re A$ and $\Im A$, respectively.

Suppose that $\Re A=\R.1$. Then $B_A(a,b):=\Re(a\bar b),\ a,b\in A$ is a \textit{standard bilinear form} on $A$ and  $Q_A(a):=B_A(a,a)$ is a \textit{standard quadratic form}. Finally, we define a multi-linear form on $\Im A$:
\begin{equation}\label{tautological form on A}
\omega_A(a,b,c):=B_A(a.b,c),\ a,b,c\in\Im A.
\end{equation}

\medskip
We will consider two products on the vector space $A\oplus A$, namely
\begin{equation}\label{CD multiplication}
(a,b)\cdot(c,d)=(ac-\bar d b,b\bar c+da)
\end{equation}
and 
\begin{equation}\label{alternative multiplication}
 (a,b)\p\cdot(c,d)=(ac+\bar d b,b\bar c+da).
\end{equation}
Then $A_2:=(A\oplus A,\cdot)$ and $\p A_2:=(A\oplus A,\p\cdot)$ are real algebras with a unit $(1,0)$.
It is clear that  
\begin{equation}\label{conjugation}
\star:A\oplus A\ra A\oplus A,\ \star{(a,b)}=(\bar a,-b)
\end{equation}
is an involution.

Then (see \cite[Section 1.13]{Y}) also $\CD(A):=(A_2,\star)$ and $\pCD(A):=(\p A_2,\star)$ are $\ast$-algebras and  $\CD(A)$  is the \textit{Cayley-Dickson algebra associated to} $(A,\ast)$. 

\bigskip

It is well known that $\C=\CD(\R),$ the algebra of quaternions $\Q=\CD(\C)$ and the algebra of octonions $\Oc=\CD(\Q)$.
\smallskip

The algebra  of pseudo-quaternions  $\pQ=\pCD(\C)$ is isomorphic to $(\mathrm M(2,\R),\ast)$ where 
\begin{equation}\label{involution on PQ}
\ast\left(
\begin{array}{cc}
a&b\\
c&d\\
\end{array}
\right)=
\left(
\begin{array}{cc}
d&-b\\
-c&a\\
\end{array}
\right).
\end{equation}
The standard quadratic form is the determinant. In particular, the signature of  $Q_{\pQ}$ is  $(2,2)$. 

We will use the following conventions. We put 
\begin{equation}\label{pseudo-quaternions}
\ \tilde I:=\left(
\begin{array}{cc}
0&-1\\
1&0\\
\end{array}
\right),\
 \tilde J:=
\left(
\begin{array}{cc}
0&1\\
1&0\\
\end{array}
\right),
\ \tilde K:=
\left(
\begin{array}{cc}
1&0\\
0&-1\\
\end{array}
\right).
\end{equation}
Be aware that $\pI.\pJ=-\pK$. 
\smallskip

The algebra of pseudo-octonions $\pOc=\CD(\pQ)\cong \pCD(\Q)$, see \cite[Section 1.13]{Y}.

\subsection{Subgroups of $\G_2$ and $\tilde\G_2$}\label{section subgroups} 
We have $\G_2=\{\varphi\in\GL(\Oc)|\ \forall a,b\in\Oc: \varphi(a.b)=\varphi(a).\varphi(b)\}$.  It is well known that  $\G_2\subset\Stab(1)\cap\Stab(\Im\Oc)$. Hence, it is more natural to view $\G_2$  as a subgroup of $\GL(\Im\Oc)$. It is well known that $\omega_\Oc\in\Lambda^3\Im\Oc^\ast$  and $\G_2=\Stab(\omega_\pOc)$. If we replace in the definition of $\G_2$ the algebra $\Oc$ by $\pOc$, then we get the group $\pG_2$. As above, $\omega_\pOc\in\Lambda^3\Im\pOc^\ast$ and that $\pG_2=\Stab(\omega_\pOc)$. We view also $\pG_2$ as a subgroup of $\GL(\Im\pOc)$.

\medskip
Let $\Sp(1)$ be the group of unit quaternions. We define
\begin{equation}\label{SO434}
\phi_{p,q}\in\GL(\Im\Q\oplus\Q),\ \phi_{p,q}(a,b):=(p.a.\bar p,q.b.\bar p),\ p,q\in\Sp(1).
\end{equation}
It is well known (see \cite[Proposition 1.20.1]{Y}) that $\SO(4)_{3,4}:=\{\phi_{p,q}:\ p,q\in\Sp(1)\}$ is a subgroup of $\G_2$. If we view $\pOc$ as $\pCD(\Q)$, then it is shown in \cite[Section 1.13]{Y} that $\SO(4)_{3,4}\subset\pG_2$. 
\medskip

We put 
\begin{equation*}
\pSp(1)_{\pm}:=\{A\in\pQ: Q_{\pQ}(A)=\pm1\}\ \ \mathrm{and}\ \ \pSp(1):=\{A\in\pQ: Q_{\pQ}(A)=1\}.
\end{equation*}
Let $A,B\in\pSp(1)_\pm$ and consider 
\begin{equation}\label{SO2234}
\tilde\phi_{A,B}\in\GL(\Im\pQ\oplus\pQ),\ \tilde\phi_{A,B}(X,Y)=Q_{\pQ}(B).(A.X.\bar A,B.Y.\bar A).
\end{equation}
It can be verified directly that $\SO(2,2)_{3,4}:=\{\tilde\phi_{A,B}|\ A,B\in\pSp(1)_\pm:\ Q_{\pQ}(A)=Q_{\pQ}(B)\}$ is a subgroup of $\pG_2$. This uses the same computation which shows that $\SO(4)_{3,4}\subset\G_2$ together with the identity $A^{-1}=Q_{\pQ}(A).\bar A,\ A\in\pSp(1)_\pm$. 

On the other hand, the group $\wt\SO(2,2)_{3,4}:=\{\tilde\phi_{A,B}|\ A,B\in\pSp(1)_\pm\}$ is not contained in $\pG_2$.  Nevertheless, we still have the following.

\begin{lemma}\label{lemma invariance of omega 5}
Let $X\in\Im\pQ,\ Y,Z\in\pQ$ and $\p\phi_{A,B}\in\wt\SO(2,2)_{3,4}$. Then
\begin{equation}
\p\phi_{A,B}.\omega_\pOc((X,0),(0,Y),(0,Z)) =\omega_\pOc((X,0),(0,Y),(0,Z))=
\end{equation}
\end{lemma}

\begin{proof}
The right hand side is
\begin{eqnarray}\label{help I}
&&\omega_\pOc((X,0),(0,Y),(0,Z))=B_\pOc((X,0).(0,Y),(0,Z))=\Re((0,Y.X).\overline{(0,Z)})\nonumber\\
&&\ \ =\Re((0,Y.X).(0,-Z))=\Re(\bar Z.Y.X,0)=\Re(\bar Z.Y.X).	
\end{eqnarray}
The left hand side is
\begin{eqnarray*}
&&\omega_\pOc(Q_{\pQ}(B)(\bar A.X.A,0),Q_{\pQ}(A)(0,\bar B.Y. A),Q_{\pQ}(A)(0,\bar B.Z.A))\\
&&=Q_{\pQ}(B)\Re(\bar A.\bar Z .B.\bar B.Y.A.\bar A.X.A)=\Re(\bar Z .Y.X).
\end{eqnarray*}
\end{proof}

The restriction map $\phi_{p,q}\mapsto\phi_{p,q}|_{\Im\Q}$ induces a homomorphism $\pi:\SO(4)_{3,4}\ra\SO(\Im\Q)$. It is well known that the homomorphism is surjective.

If $A,B\in\pSp(1)_\pm$ and $X\in\Im\pQ$, then $Q_{\pQ}(Q_{\pQ}(B) A.X.\bar A)=Q_\pQ(X)$. It follows that the restriction map $\p\phi_{A,B}\mapsto\p\phi_{A,B}|_{\Im\pQ}$ induces a map $\p\pi:\wt\SO(2,2)_{3,4}\ra\gO(\Im\pQ)$.

\begin{lemma}\label{lemma split ses of groups}
There are split short exact sequence of Lie groups
\begin{align}\label{split ses with SO432}
 &0\ra\Sp(1)\ra\SO(4)_{3,4}\xrightarrow{\pi}\SO(\Im\Q)\ra0\ \ \mathrm{and}\\
 & 0\ra\pSp(1)\ra\wt\SO(2,2)_{3,4}\xrightarrow{\tilde\pi}\gO(\Im\pQ)\ra0\label{split ses with O2234}
\end{align}
where  $\pi$ and $\tilde\pi$ are defined above. The second short exact sequence contains a split short exact sequence
\begin{equation}
0\ra\pSp(1)\ra\SO(2,2)_{3,4}\ra\SO(\Im\pQ)\ra0 \label{split SO2234}.
\end{equation}
\end{lemma}
\begin{proof}
We have $\ker\pi=\{\phi_{\pm1,q}: q\in\Sp(1)\}$. Now $\phi_{p,q}=\phi_{p',q'}$ if, and only if $p=p',q=q'$ or $p=-p',q=-q'$. It follows that $\ker\pi=\{\phi_{1,q}:\ q\in\Sp(1)\}\cong\Sp(1)$. It remains to show that the short exact sequence is split. For this consider the group $\SO(3)\cong\{\phi_{p,p}:\ p\in\Sp(1)\}$. Then  $\pi|_{\SO(3)}$ induces isomorphism $\SO(3)\ra\SO(\Im\Q)$. The inverse of this map is a splitting of $\pi$. 

It is easy to see that $\p\phi_{A,B}=\p\phi_{A',B'}$ if, and only if $A=A',B'=B$ or $A=-A',B=-B'$. We have $\ker\p\pi=\{\p\phi_{\pm1,B}:\ B\in\pSp(1)\}=\{\p\phi_{1,B}:\ B\in\pSp(1)\}\cong\pSp(1)$.  It is well known that any element of $\gO(\Im\pQ)$ is of the form $X\mapsto \pm A.X.\bar A$ for some $A\in\pSp(1)_\pm$. Then $\p\pi$ is surjective and  arguing as above, it is easy to see that it is a split. 

The last claim is easy to verify.
\end{proof}


\subsection{Maximal compact subgroup of  semi-direct product}\label{section semi-direct products}
A semi-direct product of Lie groups is  a short exact sequence 
\begin{equation}
 0\ra\gH_1\ra\gH\xrightarrow{\pi}\gH_2\ra0
\end{equation}
of Lie groups which splits, i.e. there is a homomorphism $\iota:\gH_2\ra\gH$ such that $\pi\circ\iota=Id_{\gH_2}$. In this way, we can view $\gH_2$ as a subgroup of $\gH$ and we will do this without further comment.  We write $\gH=\gH_1\rtimes\gH_2$.

\begin{lemma}\label{lemma max compact in product}
Let $\G$ be a Lie group with a compact subgroup $\gK$. Suppose that $\gH_1,\gH_2$ are closed subgroups of $\G$ such that the group generated by $\gH_1,\gH_2$ is a semi-direct product $\gH_1\rtimes\gH_2$. Assume that $\gK_i:=\gK\cap\gH_i,\ i=1,2$ is a maximal compact subgroup of $\gH_i$. 

Then the group generated by $\gK_1,\gK_2$ is a semi-direct product $\gK_1\rtimes\gK_2$ and this is a maximal compact subgroup of $\gH_1\rtimes\gH_2$.
\end{lemma}
\begin{proof}
Let $k_i\in\gK_i,\ i=1,2$ be arbitrary. Then $k_2^{-1}k_1k_2\in\gH_1\cap\gK=\gK_1$ and so the first claim follows. The group $\gK_1\rtimes\gK_2$ is homeomorphic to $\gK_1\times\gK_2$ and so it  is compact. It remains to show that it is  maximal among all compact subgroups of $\gH_1\rtimes\gH_2$.

Let $\pi:\gH_1\rtimes\gH_2\ra\gH_2$  be the canonical projection and  $\iota:\gH_2\ra\gH_1\rtimes\gH_2$ be  the inclusion. Put $\bar h:=\iota\circ\pi(h),\ h\in\gH_1\rtimes\gH_2$. Let $\gL$ be a compact group such that  $\gK_1\rtimes\gK_2\subset\gL\subset\gH_1\rtimes\gH_2$. Then in particular,  $\gK_i\subset\gL,\ i=1,2$. It is clearly enough to show that for each  $\ell\in\gL:\ \ell.\bar \ell^{-1}\in\gK_1$ and $\bar \ell\in\gK_2$.

The group  $\pi(\gL)$ is a compact subgroup of $\gH_2$ which contains $\gK_2$. By the maximality of $\gK_2$, $\pi(\gL)=\gK_2$. Thus $\bar \ell\in\iota(\pi(\gL))=\iota(\gK_2)=\gK_2$.

It is clear that $\ell.\bar \ell^{-1}\in\ker\pi=\gH_1$. As $\bar \ell\in\gK_2\subset\gL$, we have that $\ell.\bar \ell^{-1}\in\gL$ and thus $\ell.\bar \ell^{-1}\in\gL\cap\gH_1$. Now $\gL\cap\gH_1$ is a compact subgroup of $\gH_1$ which contains $\gK_1$. Again by the maximality of $\gK_1$, $\gK_1=\gL\cap\gH_1$ and  the proof is complete.
\end{proof}

\begin{lemma}\label{lemma max compact in product I}
Let $\gH_1\rtimes\gH_2$ be a semi-direct product of Lie groups. Assume that the trivial group $\{e\}$ is a maximal compact subgroup of $\gH_i$  and that $\gK$ is a maximal compact subgroup of $\gH_j$ where $i,j\in\{1,2\},\ i\ne j$. Then $\gK$ is a maximal  compact subgroup of $\gH_1\rtimes\gH_2$.
\end{lemma}
\begin{proof}
As $\gH_1\rtimes\gH_2$ is diffeomorphic to $\gH_1\times\gH_2$, the groups $\gH_i,\ i=1,2$ are closed in the semi-direct product. The group $\gK$ is a compact subgroup of $\gH_1\rtimes\gH_2$. By the assumptions,  $\gK\cap\gH_i$ is a maximal compact subgroup of $\gH_i,\ i=1,2$.  The claim follows from Lemma \ref{lemma max compact in product}
\end{proof}

\begin{lemma}\label{lemma max compact in upper trian}
Let $\gL$ be a subgroup of $\GL(n,\R)$. Suppose that  for each $A\in\gL$ the matrix $A-1_n$ is  strictly upper triangular. Then  $\{1_n\}$ is a maximal compact subgroup of $\gL$.
\end{lemma}
\begin{proof}
Let $\gK$ be a maximal compact subgroup of $\gL$ and $B$ be a $\gK$-invariant inner product on $\R^n$. Let $\{f_1,\dots,f_n\}$ be a $B$-orthonormal basis that is obtained by applying the Gram-Schmidt algorithm to the standard basis of $\R^n$. Let $F=(f_1|f_2|\dots|f_n)$ be the corresponding matrix. Then $F^{-1}.\gK.F\subset\gO(7)$. On the other hand,   $F$ and $F^{-1}$ are upper triangular matrices with positive numbers on the diagonal. Thus, for each $A\in\gK$ the matrix $F^{-1}.A.F$  is  upper triangular and its diagonal coefficients are equal to 1. Hence,  $F^{-1}.A.F=1_n$ and  $A=1_n$.  
\end{proof}

\section{Maximal compact subgroups}\label{section max compact subgroups}
This Section is organized as follows. Each type is treated in its own section. We first recall some results about $\gO_i:=\Stab(\omega_i)$ from \cite{BV} where  $\omega_i$ is a fixed  representative  of the type $i=1,\dots,8$. Then we use these results to find a maximal compact subgroup of $\gO_i$. We separate the summary part from the genuine research part of this article by \textemdash\ x \textemdash.
\medskip

\subsection{Type 8}\label{section type 8}
A representative chosen in \cite[Section Type 8.]{BV} is
\begin{equation}\label{omega8}
\omega_8=\alpha_{123}+\alpha_{145}-\alpha_{167}+\alpha_{246}+\alpha_{257}+\alpha_{347}-\alpha_{356}.
\end{equation}
It is well known that  $\gO_8=\G_2$.

\begin{center}
 \textemdash\ x \textemdash
\end{center}

We will need in Section \ref{section type 7} the following observation.
Recall (see Section \ref{section algebras}) that  $\Oc=\CD(\Q)$ and  the definition (\ref{tautological form on A}) of the tautological 3-form $\omega_\Oc$.
Let $\Phi_8:\mV\ra\Im\Oc$  be the linear isomorphism that maps the standard basis of $\mV$ to the standard basis 
\begin{equation}
\label{standard basis of O}
\{(i,0),(j,0),(k,0),(0,1),(0,i),(0,j),(0,k)\}. 
\end{equation}
of $\Im\Oc$ where $\{1,i,j,k\}$ is the standard basis of $\Q$. 

\begin{lemma}\label{lemma type 8}
$\omega_8=\Phi_8^\ast\omega_\Oc$.
\end{lemma}
\begin{proof}
This is a straightforward computation which uses (\ref{help I}) together with $\omega_\Oc((i,0),(j,0),(k,0))=\Re((k,0).\overline{(k,0)})=\Re(-k^2)=\Re(1)=1.$ 
\end{proof}
\medskip

\subsection{Type 7}\label{section type 7}
Let us first recall some facts from \cite[Section Type 7.]{BV}.
A representative is 
\begin{equation}
\omega_7=\alpha_{125}+\alpha_{136}+\alpha_{147}+\alpha_{237}-\alpha_{246}+\alpha_{346}.
\end{equation}
Then $\mV_3:=\triangle_7^3=[e_5,e_6,e_7]$ and  we put $\mW_4:=\mV/\mV_3$. The 3-form $\omega_7$ induces a bilinear form\footnote{Here $v=(c_1,\dots,c_7),\ w=(d_1,\dots,d_7)$} $B(v,w)=c_1d_1+c_2d_2+c_3d_3+c_4d_4$ on $\mV$.  The form $B$ descends to a positive definite bilinear form $B_4$ on $\mW_4$.  The associated quadratic form is denoted by $Q_4$. We have\footnote{Recall notation from (\ref{notation stabilizer of subspace}). In particular, $\Stab([B])$ is the stabilizer of the \textbf{line} spanned by $B$.} $\gO_7\subset\Stab(\mV_3)\cap\Stab([B])$. 

The insertion map $v\mapsto i_v\omega_7$ induces a linear map $\lambda:\mV_3\ra\Lambda^2\mW_4$. We put $\sigma_i:=\lambda(e_{i+4}),\ i=1,2,3$. As $B_4$ is non-degenerate, there are unique $E,F,G\in\End(\mW_4)$ such that 
\begin{equation}\label{type 7 quaternionic triple}                                                                                                                                                                                                                            
\sigma_1=B_4\circ E,\ \sigma_2=B_4\circ F,\  \sigma_3=B_4\circ G
\end{equation}
where $(B_4\circ A)(u,v):=B_4(A(u),v),\ A\in\End(\mW_4),\ u,v\in\mW_4$. The endomorphisms satisfy $E^2= F^2=-Id_{W_4}$ and $E F=-F E=G$. Hence, $\mW_4$ is a 1-dimensional  $\Q$-vector space. 

The scalar product $B_4$ induces a scalar product on $\Lambda^2\mW_4^\ast$. As $\lambda$ is injective, $B_4$ induces a scalar product $B_3$ on $\mV_3$ and one finds that $\{e_5,e_6,e_7\}$ is an orthonormal basis.

The homomorphism $\rho:\gO_7\ra\GL(\mV_3),\ \varphi\mapsto \varphi|_{\mV_3}$ induces a split short exact sequence $0\ra\gK\ra\gO_7\ra\CSO(Q_3)\ra0$. Any $\varphi\in\gO_7$ induces an endomorphism $\p\varphi$ of $\mW_4$ which preserves $Q_4$ up to scale. If $\varphi\in\gK$, then $\tilde\varphi\in\GL_\Q(\mW_4)$ and  $\tilde\varphi\in\Stab(B_4)$. This induces another split short exact sequence  $0\ra\gL\ra\gK\ra\Sp(1)\ra0$. The group $\gL$ is isomorphic to the Lie group $\R^8$. Notice that $ 
\ell_{\cB_{st}}=\left(
\begin{matrix}
1_4&0\\
\ast&1_3\\
\end{matrix}
\right)$ for every $\ell\in\gL$.
All together, $\gO_7$ is isomorphic to a  semi-direct product 
\begin{equation}\label{semi-direct product type 7}
 (\R^8\rtimes\Sp(1))\rtimes\CSO(3)
\end{equation}
as claimed in \cite[10. Proposition]{BV}

\begin{center}
 \textemdash\ x \textemdash
\end{center}

We can now proceed. We put $$f_1:=-e_6,\ f_2:=e_7,\ f_3:=e_5, \ f_4:=-e_1,\ f_5:=e_3,\ f_6:=-e_4,\ f_7:=-e_2$$
and  $\cB_7:=\{f_1,\dots,f_7\}$.  Let $\cB_7^\ast=\{\beta_1,\dots,\beta_7\}$ be a dual basis. Then
\begin{equation}
\omega_7=\beta_{145}-\beta_{167}+\beta_{246}+\beta_{257}+\beta_{347}-\beta_{356}.
\end{equation}

Consider a linear isomorphism $\Phi_7:\mV\ra Im\Oc$ which maps the basis $\cB_7$  to the standard basis (\ref{standard basis of O}) of $\Im\Oc$.  Recall Section \ref{section subgroups} that  $\SO(4)_{3,4}$ is a subgroup of $\GL(\Im\Oc)$. Thus, $\gK_7:=\Phi_7^\ast\SO(4)_{3,4}$ is a subgroup of $\GL(\mV)$. Here we use the notation set in (\ref{notation pullback of map}).

\begin{lemma}\label{lemma inclusion type 7}
The subspaces $\mV_3$ and $\mV_4:=[f_4,\dots,f_7]$ are $\gK_7$-invariant and 
$\gK_7\subset\gO_7$.
\end{lemma}
\begin{proof}
Since $\mV_3=\Phi_7^{-1}(\Im\Q)$ and $\mV_4=\Phi_7^{-1}(\Q)$, the first claim follows from the definition (\ref{SO434}).

By Lemma \ref{lemma type 8}: $\omega_7=\Phi_7^\ast\omega_\Oc-\beta_{123}$. As $\SO(4)_{3,4}\subset\G_2$, we see that $\gK_7\subset\Phi_7^\ast\G_2=\Stab(\Phi_7^\ast\omega_\Oc)$. From (\ref{SO434}) is also clear that $\gK_7\subset\SO(B_3)\subset\Stab(\beta_{123})$. Thus $\gK_7\subset\gO_7$.
\end{proof}

Put
\begin{equation}\label{group R+}
\gR^+:=\bigg\{
\varphi\in\GL(\mV)\bigg|
\  \varphi_{\cB_7}= 
\left(
\begin{matrix}
\lambda^{-2}.1_3&0\\
0&\lambda.1_4\\
 \end{matrix}
\right),\lambda>0
\bigg\}
\end{equation}
It is clear that $\gR^+$ is a subgroup of $\GL(\mV)$ and that  $\gR^+\subset\gO_7$. From Lemma \ref{lemma inclusion type 7} follows that $\gK_7$ commutes with $\gR^+$. It is  easy to see that $k.\ell.k^{-1}\in\gL$ whenever $\ell\in\gL,\ k\in\gK_7\times\gR^+$. Hence, the group generated by $\gL,\gK_7$ and $\gR^+$ is isomorphic to $\gL\rtimes(\gK_7\times\gR^+)$. This is a subgroup of $\gO_7$.

\begin{lemma}\label{lemma semi-direct product type 7}
$\gO_7=\gL\rtimes(\gK_7\times\gR^+)$. 
\end{lemma}
\begin{proof}
Let $v\in\mV_3$. The map $\Phi_7$ restricts to an isomorphism $\mV_3\ra\Im\Q$ and so there is a unique $a\in\Im\Q$ so that  $\Phi_7(v)=(a,0)$.
The composition $\mV\xrightarrow{\Phi_7}\Im\Oc\ra\Q$ where the second map is the canonical projection descends to an isomorphism $\underline\Phi_7:\mW_4\ra\Q$. Notice that $\underline\Phi_7^\ast B_\Q=B_4$. We know that there is a unique $A\in\End(\mW_4)$ such that $\lambda(v)=B_4\circ A$. We will now show that $A$ corresponds to the multiplication by  $a$ on the right.

Let $w_1,w_2\in\mW_4$ and put  $\ b=\underline\Phi_7(w_1),\ c=\underline\Phi_7(w_2)$.  Then
\begin{eqnarray*}
&&\lambda(v)(w_1,w_2)=\omega_\Oc((a,0),(0,b),(0,c))=\Re(\bar cba)=\Re(ba \bar c)=B_\Q(b.a,c)\\
&&=B_4(\Phi_7^\ast(r_a)(w_1),w_2)
\end{eqnarray*}
where $r_a:\Q\ra\Q$ is  the right-action by $a$. We see that $A=\underline\Phi_7^\ast(r_a)$. 

Recall Lemma \ref{lemma split ses of groups} that $\pi:\SO(4)_{3,4}\ra\SO(\Im\Q),\ \varphi\mapsto\varphi|_{\Im\Q}$ induces the split short exact sequence (\ref{split ses with SO432}) where $\Sp(1)=\{\phi_{1,q}:\ q\in\Sp(1)\}$ and that $\phi_{1,q}(a,b)=(a,q.b)$. Now it is easy to see that $\gK_7\times\gR^+$ is a splitting of the subgroup $\Sp(1)\rtimes\CSO(3)$ from (\ref{semi-direct product type 7}). The claim then follows from the definition of $\gL$.
\end{proof}

\begin{thm} \label{thm max compact type 7}
$\gK_7$ is a maximal compact subgroup of $\gO_7$.
\end{thm}
\begin{proof}
By Lemma \ref{lemma max compact in upper trian},  $\{1_7\}$ is a maximal compact subgroup of $\gL$. It is also clear that $\{1_7\}$ is a maximal compact subgroup of $\gR^+$. The claim follows from Lemma \ref{lemma max compact in product I}.
\end{proof}

\medskip

\subsection{Type 5}\label{section type 5}
A representative chosen in \cite[Section type 5]{BV} is
 \begin{eqnarray}\label{type 5}
\omega_5=\alpha_{123}-\alpha_{145}+\alpha_{167}+\alpha_{246}+\alpha_{257}+\alpha_{347}-\alpha_{356}.
\end{eqnarray}
It is well known that $\gO_5=\pG_2$.
\begin{center}
\textemdash\ x \textemdash
\end{center}

We will later need two more representatives of type 5.
Recall Section \ref{section algebras} that  $\pOc=\pCD(\Q)$ and  the definition of  $\omega_\pOc$. Let $\Phi_5:\mV\ra\Im\pOc$ be the linear isomorphism that maps the standard basis of $\mV$ to the  basis (\ref{standard basis of O}) of $\Im\pOc$. From  (\ref{help I}) and (\ref{alternative multiplication})  follows at once that $\Phi_5^\ast\omega_\pOc=2\alpha_{123}-\omega_8$, i.e.
 \begin{eqnarray}\label{type 5.I}
\Phi^\ast_5\omega_\pOc=\alpha_{123}-\alpha_{145}+\alpha_{167}-\alpha_{246}-\alpha_{257}-\alpha_{347}+\alpha_{356}.
\end{eqnarray}
Notice that  $\varphi.\omega_5=\Phi^\ast_5\omega_\pOc$ where $\varphi\in\GL(\mV)$ is determined by $\varphi(e_i)=e_i, \ 1\le i\le 5$ and $\varphi(e_j)=-e_j,\ j =6,7$.

\bigskip

Let us now view $\pOc$ as $\CD(\pQ)$ and $\pQ$ as $M(2,\R)$. Consider the linear isomorphism $\tilde\Phi_5:\mV\ra\Im\pOc$ that sends the standard basis of $\mV$ to the basis 
\begin{equation}\label{basis of pO}
\{(\pI,0),(\pJ,0),(\pK,0),(0,\sqrt2E_{11}),(0,\sqrt2E_{21}),(0,\sqrt2E_{12}),(0,\sqrt2E_{22})\}.
\end{equation}

\begin{lemma}\label{lemma type 5 II}
\begin{eqnarray}
\tilde\Phi_5^\ast\omega_\pOc=\alpha_{123}+\alpha_{145}+\alpha_{167}-\alpha_{245}+\alpha_{267}+\alpha_{347}-\alpha_{356}.
\end{eqnarray}
\end{lemma}
\begin{proof}
Let $\tr$ be the usual trace and
\begin{equation}
A=\left(
\begin{matrix}
a_3&a_2-a_1\\
a_1+a_2&-a_3
\end{matrix}
\right), \
B=\sqrt2\left(
\begin{matrix}
b_1&b_3\\
b_2&b_4
\end{matrix}
\right),\ C=\sqrt2\left(
\begin{matrix}
c_1&c_3\\
c_2&c_4
\end{matrix}
\right).
\end{equation} 
From  (\ref{help I}) follows that
\begin{align*}
\omega_\pOc((A,0),(0,B),(0,C))&=\Re(\bar CBA)=\frac{1}{2}\tr(\bar CBA)\\
&=a_1(b_1c_2-b_2c_1+b_3c_4-b_4c_3)\\
&+a_2(-b_1c_2+b_2c_1+b_3c_4-b_4c_3)\\
&+a_3(b_1c_4-b_4c_1-b_2c_3+b_3c_2).
\end{align*}
As $\omega_\pOc((\pI,0),(\pJ,0),(\pK,0))=\Re((-\pK,0).\overline{(\pK,0)})=\Re(\pK^2)=\Re(1)=1,$
the claim follows.
 \end{proof}

\medskip

\subsection{Type 2}\label{section type 2}
Let us first recapitulate  \cite[Section Type 2.]{BV}.
A representative   is
\begin{equation*}
\omega_2=\alpha_{125}+\alpha_{127}+\alpha_{147}-\alpha_{237}+\alpha_{346}+\alpha_{347}.
\end{equation*}
We will use a  basis $\cB_2=\{f_5,f_6,f_7,e_1,e_2,e_3,e_4\}$ of $\mV$ where
$$f_5:=e_5+e_6,\ f_6:=-e_5+e_6,\ f_7:=-e_5-e_6+e_7.$$
For a dual basis $\cB_2^\ast=\{\beta_1,\ldots,\beta_7\}$ we have 
$\beta_1=\frac{1}{2}(\alpha_5+\alpha_6+2\alpha_7),\ \beta_2=\frac{1}{2}(-\alpha_5+\alpha_6),\ \beta_3=\alpha_7$ and $\beta_{i+3}=\alpha_i,\ i=1,\dots,4$
and  so
$$\alpha_i=\beta_{i+3},\ i=1,\dots,4,\ \alpha_5=\beta_1-\beta_2-\beta_3,\ \alpha_6=\beta_1+\beta_2-\beta_3,\ \alpha_7=\beta_3.$$
Then 
\begin{equation}
\omega_2=\beta_{145}+\beta_{167}-\beta_{245}+\beta_{267}+\beta_{347}-\beta_{356}.
\end{equation}

The subspace  $\mV_3:=[f_5,f_6,f_7]$ is $\gO_2$-invariant and we put $\mW_4:=\mV/\mV_3$. We denote by $B_3$ a symmetric bilinear form on $\mV_3$ with orthogonal basis $\{f_5,f_6,f_7\}$ and $B_3(f_5,f_5)=1,\ B_3(f_6,f_6)=B_3(f_7,f_7)=-1$. The 3-form induces a  quadratic form $Q(v)=2(c_1c_4-c_2c_3)$  where $v=(c_1,\dots,c_7)$. The associated symmetric bilinear form $B$ descends to a bilinear form $B_4$ on $\mW_4$. We denote by $Q_3$ and $Q_4$ the quadratic form on $\mV_3$ and $\mW_4$ associated to $B_3$ and $B_4$, respectively. We have $\gO_2\subset\Stab([B_3])\cap\Stab([B])$.

The  insertion map $v\mapsto i_v\omega_2$ induces an injective linear map $\lambda:\mV_3\ra\Lambda^2\mW_4$. We put $\sigma_i:=\lambda(f_{4+i}),\ i=1,2,3$. As $B_4$ is non-degenerate, there are unique  $E,F,G\in\End(\mW_4)$ such that 
\begin{eqnarray*}                                                                                                                                                                                                                            &&\sigma_1= B_4\circ E,\ \sigma_2=B_4\circ F,\ \sigma_3=B_4\circ G                                                                                                                                                                                                          \end{eqnarray*}
where we use the notation from (\ref{type 7 quaternionic triple}).
Then $E^2=-F^2=-Id_{\mW_4}$ and $EF=-FE=G$ and hence, $\mW_4$ is a 1-dimensional free $\pQ$-module.

There is a split short exact sequence $0\ra\gO_2^+\ra\gO_2\xrightarrow{\det}\R^\ast\ra0$. The restriction map $\varphi\in\gO^+_2\ra\varphi|_{\mV_3}$ induces another split short exact sequence $0\ra\gK\ra\gO^+_2\ra\SO(Q_3)\ra0$. Next, any $\varphi\in\gO_2$ induces $\p\varphi\in\GL(\mW_4)$. If $\varphi\in\gK$, then $\p\varphi\in\GL_\pQ(\mW_4)\cap\Stab(B_4)\cong\widetilde\Sp(1)$. We obtain another split short exact sequence $0\ra\gL\ra\gK\ra\pSp(1)\ra0$. The group $\gL$ is isomorphic to $\R^8$. Notice that   
$\ell_{\cB_2}=\left(
 \begin{matrix}
1_3&\ast\\
0&1_4\\
 \end{matrix}
 \right)$
for every $\ell\in\gL$. All together,  $\gO_2$ is isomorphic to a semi-direct product
\begin{equation}\label{semi-direct product type 2}
((\R^8\rtimes\pSp(1))\rtimes\SO(1,2))\rtimes\R^\ast
\end{equation}
as claimed in \cite[5. Proposition]{BV}.

\begin{center}
 \textemdash\ x \textemdash
\end{center}

We can now continue. Let $\Phi_2:\mV\ra\Im\pOc$ be a linear isomorphism that sends  $\cB_2$ to the basis (\ref{basis of pO}). Recall  Section \ref{section subgroups} that $\SO(2,2)_{3,4}\subset\wt\SO(2,2)_{3,4}\subset\GL(\Im\pOc)$. Put $\wt\gH_2:=\Phi_2^\ast\wt\SO(2,2)_{3,4}$ and $\gH_2:=\Phi_2^\ast\SO(2,2)_{3,4}$.

\begin{lemma}\label{lemma inclusion type 2}
The subspaces $\mV_3$ and $\mV_4:=[e_1,\dots,e_4]$ are $\wt\gH_2$-invariant and  
$\gH_2\subset\gO_2$. 
\end{lemma}
\begin{proof}
 The first claim follows from the fact that $\Phi_2(\mV_3)=\Im\pQ,\ \Phi_2(\mV_4)=\pQ$  are $\wt\SO(2,2)_{3,4}$-invariant subspaces. 
From Lemma \ref{lemma type 5 II} follows  that $\omega_2=\Phi_2^\ast\omega_{\pOc}-\beta_1\wedge\beta_2\wedge\beta_3$. The second claim is then a consequence of Lemma \ref{lemma invariance of omega 5}.
\end{proof}

Let $\gR^+$ be the group from (\ref{group R+}). It is easy to see that $\gR^+\subset\gO_2$. From Lemma \ref{lemma inclusion type 2} follows that  $\gR^+$ commutes with $\wt\gH_2$. Using the definition of $\gL$, one can easily check that $g.\ell.g^{-1}\in\gL$ whenever $\ell\in\gL,\ g\in\wt\gH_2\times\gR^+$. We see that $\gL,\wt\gH_2$ and $\gR^+$  generate a group $\gL\rtimes(\wt\gH_2\times\gR^+)$ inside $\gO_2$.

\begin{lemma}\label{lemma semi-direct product type 2}
$\gO_2=\gL\rtimes(\wt\gH_2\times\gR^+)$.
\end{lemma}
\begin{proof}
Let $v\in\mV_3$. The map $\Phi_2$ restricts to an isomorphism $\mV_3\ra\Im\pQ$ and so there is a unique $a\in\Im\pQ$ so that  $\Phi_2(v)=(a,0)$.
The composition $\mV\xrightarrow{\Phi_2}\Im\pOc\ra\pQ$ where the second map is the canonical projection descends to an isomorphism $\underline\Phi_2:\mW_4\ra\pQ$.  From the summary given above, there is a unique $A\in\End(\mW_4)$ such that $\lambda(v)=B_4\circ A$.
Following the proof of Lemma \ref{lemma semi-direct product type 7}, we find that  $A=\underline\Phi_2^\ast(r_a)$.  

It is easy to see that $\det:\gL\rtimes(\wt\gH_2\times\gR^+)\ra\R^\ast$ is surjective and that $\ker(\det)=\gL\rtimes\gH_2$.
Recall Lemma \ref{lemma split ses of groups} that  $\p\pi$ induces the split short exact sequence (\ref{split ses with SO432}) where  $\pSp(1)=\{\p\phi_{1,B}:\ B\in\pSp(1)\}$ and  $\p\phi_{1,B}(X,Y)=(X,B.Y)$. From this easily follows that $\gH_2$ is a splitting of  $\pSp(1)\rtimes\SO(1,2)$ in (\ref{semi-direct product type 2}). The claim  now follows from the definition of $\gL$.
\end{proof}

It remains to find a maximal compact subgroup of $\gO_2$. We will for simplicity consider only a maximal connected and compact subgroup.
The group $\SO(2)$ is a  maximal compact  subgroup of $\SL(2,\R)\cong\pSp(1)$. From the proof of Lemma \ref{lemma split ses of groups}  follows that a subgroup $\gK^0$ of $\SO(2,2)_{3,4}$ that is generated by  $\SO(2)_1:=\{\p\phi_{1,B}:\ B\in\SO(2)\}$ and $\SO(2)_2:=\{\p\phi_{A,A}:\ A\in\SO(2)\}$ is a semi-direct product $\gK^o:=\SO(2)_1\rtimes\SO(2)_2$. As $\SO(2)$ is commutative, the product is direct. Put $\gK^o_2:=\Phi_2^\ast\gK^o$.

\begin{thm}\label{thm max compact type 2}
The group $\gK^o_2\cong\SO(2)\times\SO(2)$ is a maximal connected and compact subgroup of $\gO_2$. There is an isomorphism of $\gK_2^o$-modules
\begin{equation}
\mV\cong\R\oplus\C_1\oplus\C_2\oplus\C_1\otimes_\C\C_2 
\end{equation}
where we denote by $\C_i,\ i=1,2$  the standard  complex representation of the $i$-th factor of $\SO(2)\times\SO(2)$ on $\C=\R^2$.
\end{thm}
\begin{proof}
By Lemma \ref{lemma max compact in upper trian}, $\{1_7\}$ is a maximal compact subgroup of $\gL$. It is also clear that $\{1_7\}$ is a maximal compact subgroup of $\gR^+$. Hence, a maximal compact subgroup of $\wt\gH_2$ is a maximal compact subgroup of $\gO_2$. 
From Lemma \ref{lemma max compact in product} easily follows that $\gK^o$ is a maximal compact  subgroup of the connected component of the identity element of $\SO(2,2)_{3,4}$. 
Hence, $\gK^o_2$ is a maximal compact subgroup  of the connected component $\gO_2^o$ of the identity element  of $\gO_2$. By the Cartan-Iwasawa-Malcev theorem, any two  maximal compact  subgroups of $\gO_2^o$ are conjugated and hence isomorphic. The first claim follows and it remains to show the second claim.

 Put 
$$\cB_{\pOc}:=\{(\pI,0),(\pJ,0),(-\pK,0),(0,1),(0,\pI),(0,\pJ),(0,-\pK)\}.$$ 
Let $R_t:\R^2\ra\R^2$ be the anti-clockwise rotation at angle $t\in\R$. Then it is straightforward to find that
\begin{equation}\label{element of max compact type 2}
(\p\phi_{1,R_s})_{\cB_{\pOc}}=
\left(
 \begin{matrix}
1&0&0&0\\
0&1&0&0\\
0&0&R_{s}&0\\
0&0&0&R_{s}\\
 \end{matrix}
 \right),\ (\p\phi_{R_t,R_t})_{\cB_{\pOc}}=
\left(
 \begin{matrix}
1&0&0&0\\
0&R_{2t}&0&0\\
0&0&1&0\\
0&0&0&R_{2t}\\
 \end{matrix}
 \right)
\end{equation}
and  the second claim now easily follows.
\end{proof}

\subsection{Type 6}\label{section type 6}
A representative is
\begin{equation}
 \omega_6=\alpha_{127}-\alpha_{136}+ \alpha_{145}+\alpha_{235}+\alpha_{246}.
\end{equation}
Invariant subspaces are $\triangle_6^2=[e_7]$ and $\triangle_6^3=[e_3,e_4,\dots,e_7]$. We denote them by $\mV_1$ and $\mV_5$, respectively. We put $\mW_2:=\mV/\mV_5,\ \mZ_4=\mV_5/\mV_1$. The 3-form $\omega_6$ induces a quadratic form $Q(v)=c_1^2+c_2^2,\ v=(c_1,\dots,c_7)$. Let $B$ be the associated bilinear form. Then $\gO_6\subset\Stab([B])$. The form $Q$ descends to a regular quadratic form $Q_2$ on $\mW_2$. We denote by $B_2$ the associated bilinear form. 

The insertion map $v\mapsto i_v\omega_6|_{V_5}$ induces a monomorphism $\lambda:\mW_2\ra\Lambda^2\mZ_4^\ast$ and we put $\sigma_i:=\lambda(e_i),\ i=1,2$. 

Each $\varphi\in\gO_6$ induces $\tilde\varphi\in\GL(\mW_2)$. Since $\tilde\varphi\in\Stab([B_2])$, we get a map $\mu:\gO_6\ra\CO(Q_2)$ and a split short exact sequence $0\ra\gO_6^1\ra\gO_6\xrightarrow{\det\mu}\R^\ast\ra0$. Then $\mu|_{\gO_6^1}$ induces another split short exact sequence $0\ra\gL\ra\gO_6^1\ra\SO(Q_2)\ra0$. Each $\varphi\in\gO_6$ descends to  $\bar\varphi\in\GL(\mZ_4)$ and we get a map $\rho:\gO_6\ra\GL(\mZ_4)$. If $\varphi\in\gL$, then $\bar\varphi.\sigma_i=\sigma_i,\ i=1,2$. It can be shown that  $\Stab(\sigma_1)\cap\Stab(\sigma_2)\cong\SL(2,\C)$ and we get a split short exact sequence $0\ra\gM\ra\gL\ra\SL(2,\C)\ra0$. Moreover, it can be shown that 
$\varphi_{\cB_{st}}=\left(
\begin{matrix}
1_2&0&0\\
\ast&1_4&0\\
\ast&\ast&1\\
 \end{matrix}
\right)$ for every $\varphi\in\gM$ and that  $\gM\cong(\R^6\rtimes\R^2)\rtimes\R^2$. All together, $\gO_6$ is isomorphic to a semi-direct product
\begin{equation}\label{semidirect product type 6}
((((\R^6\rtimes\R^2)\rtimes\R^2)\rtimes\SL(2,\C))\rtimes\SO(2))\rtimes\R^\ast 
\end{equation}
as claimed in \cite[9. Proposition]{BV}.

\begin{center}
 \textemdash\ x \textemdash
\end{center}

Let us first relabel the standard basis. We put
$$f_1:=e_7,\ f_2:=e_1,\ f_3:=e_2,\ f_4:=e_4,\ f_5:=-e_3,\ f_6:=-e_5,\ f_7:=-e_6.$$ Let  $\cB_6^\ast=\{\beta_1,\dots,\beta_7\}$ be a dual basis to $\cB_6:=\{f_1,f_2,\dots,f_7\}$. Then we have $$\beta_1=\alpha_7,\ \beta_2=\alpha_1,\ \beta_3=\alpha_2,\ \beta_4=\alpha_4,\ \beta_5=-\alpha_3,\ \beta_6=-\alpha_5,\ \beta_7=-\alpha_6.$$ From this one finds  that 
\begin{equation}
\omega_6=\beta_{123}-\beta_{246}-\beta_{257}-\beta_{347}+\beta_{356}.
\end{equation}

\bigskip

Consider the linear isomorphism $\C^2\ra\Q,\ (x_1+iy_1, x_2+iy_2)\mapsto(x_1+iy_1+y_2j+x_2k)$. This isomorphism is complex linear with respect to the multiplication  by $i$ on $\Q$ on the right. This induces embedding $\SL(2,\C)\hookrightarrow\GL(\Q)$. We compose this with the canonical map $\GL(\Q)\hookrightarrow\GL(\Im\Q\oplus\Q),\ \varphi\mapsto Id_{\Im\Q}\oplus\varphi$ so that we can view $\SL(2,\C)$ as a subgroup of $\GL(\Im\Q\oplus\Q)$.

Put $\widetilde\gU(1):=\{p\in\Sp(1):p.i=\pm i.p\}$. It is easy to see that $\widetilde\gU(1)$ has two connected components $\widetilde\gU(1)_0=\{e^{it}:\ t\in\R\}$ and $\widetilde\gU(1)_1=\{e^{it}.j:\ t\in\R\}$. We put $\gK_6':=\{\phi_{p,q}:\ p\in\widetilde\gU(1),\ q\in\Sp(1)\}$. It is easy to see that $\gK_6'$ is a subgroup of $\GL(\Im\Q\oplus\Q)$ and we denote by $\gH_6'$ the group generated by $\gK_6'$ and $\SL(2,\C)$.

\begin{lemma}\label{lemma split ses with H_6'}
The map $\pi_6:\gH_6'\ra\SO(\Im\Q),\ \pi_6(\varphi)=\varphi|_{\Im\Q}$ induces a split short exact sequence
\begin{equation}\label{ses with H_6'}
0\ra\SL(2,\C)\ra\gH_6'\ra\gO(2)\ra0
\end{equation}
where $\gO(2)\subset\SO(\Im\Q)$ is the stabilizer of the line spanned by $i$. The group $\gK_6'$ is a maximal compact subgroup of $\gH'_6$. 
\end{lemma}
\begin{proof}
It is easy to see that the image of $\pi_6$ is  $\gO(2)$. From Lemma \ref{lemma split ses of groups} follows that $\ker\pi_6$ is generated by $\SL(2,\C)$ and by $\Sp(1)=\{\phi_{1,q}:\ q\in\Sp(1)\}$. From (\ref{SO434}) is clear that $\Sp(1)\subset\SL(2,\C)$ and hence, $\pi_6$ induces the short exact sequence (\ref{ses with H_6'}). This sequence is split by repeating the proof of Lemma \ref{lemma split ses of groups}. It remains to prove the second claim.

The group $\gK_6'$ is a compact subgroup of $\gH_6'=\SL(2,\C)\rtimes\gO(2)$. We have that $\gK_6'\cap\gO(2)=\gO(2)$ and that $\gK_6'\cap\SL(2,\C)=\Sp(1)$ is a maximal compact subgroup of $\SL(2,\C)$. The claim  follows from  Lemma \ref{lemma max compact in product}.
\end{proof}

Let us view $\pOc$ as $\pCD(\Q)$. Let $\Phi_6:\mV\ra\Im\Q\oplus\Q$ be the linear isomorphism which maps the basis $\cB_6$ to the basis (\ref{standard basis of O}) of $\Im\pOc$. We put $\gH_6:=\Phi_6^\ast\gH_6'$  and $\gK_6:=\Phi_6^\ast\gK_6'$. These are by definition subgroups  of $\GL(\mV)$. 

\begin{lemma}\label{lemma inclusion type 6}
The subspaces $[f_1],\ [f_2,f_3],\ [f_4,\dots,f_7]$ are $\gH_6$-invariant and
$\gH_6\subset\gO_6$.
\end{lemma}
\begin{proof}
Notice that $[f_1]=\Phi_6^{-1}([(i,0)]),\ [f_2,f_3]=\Phi_6^{-1}([(j,0),(k,0)])$ and $[f_4,\dots,f_7]=\Phi_6^{-1}(\pQ)$. Hence, the first claim follows from the definition of $\gH_6'$ and  Lemma \ref{lemma split ses with H_6'}. Notice that $\gH_6\subset\Stab(\beta_1)\cap\Stab(\beta_{23})$.

The 3-form $\Phi_6^\ast\omega_{\pOc}$ is obtained from  $\Phi_5^\ast\omega_{\pOc}$ given in (\ref{type 5.I}) by replacing each $\alpha_i$ by $\beta_i,\ i=1,2,\dots,7$. We see that $\omega_6+\beta_1\wedge\gamma_1=\Phi_6^\ast\omega_{\pOc}$ where  $\gamma_1:=-\beta_4\wedge\beta_5+\beta_6\wedge\beta_7$. 

The group $\gH_6$ is generated by  $\gK_6$ and $\Phi_6^\ast\SL(2,\C)$.  Hence, it is enough to show that  $\gK_6\subset\Stab(\omega_6)$ and $\SL(2,\C)\subset\Stab(\omega_6)$.

Let  $\varphi\in\gK_6$. As $\gK_6'\subset\SO(4)_{3,4}\subset\tilde\G_2$, it follows that $\gK_6\subset\Phi_6^\ast\pG_2=\Stab(\Phi^\ast_6\omega_{\pOc})$.  Hence, it is enough to show $\varphi.\gamma_1=\gamma_1$. We have
\begin{eqnarray*}
\varphi.(\Phi_6^\ast\omega_{\pOc})&=&\varphi.(\beta_{123}+\beta_1\wedge\gamma_1-\beta_{246}-\beta_{257}+\beta_{356}-\beta_{347})\\
&=&\beta_{123}+\beta_1\wedge\varphi.\gamma_1+\varphi.(-\beta_{246}-\beta_{257}+\beta_{356}-\beta_{347})\\
&=&\beta_{123}+\beta_1\wedge\varphi.\gamma_1+\beta_2\wedge\gamma_2+\beta_3\wedge\gamma_3
\end{eqnarray*}
for some 2-forms $\gamma_2,\gamma_3$ which  belong to the subalgebra of $\Lambda^\bullet\mV^\ast$ generated $\beta_4,\beta_5,\beta_6,\beta_7$. But since $\varphi.(\Phi_6^\ast\omega_{\pOc})=\Phi_6^\ast\omega_{\pOc}$, it follows that $\varphi.\gamma_1=\gamma_1$ and $\gamma_2=-\beta_{46}-\beta_{57},\ \gamma_3=\beta_{56}-\beta_{47}$. 

As $\SL(2,\C)$ acts by identity on $\Im\Q$, we know that $\beta_i,\ i=1,2,3$ are invariant under $\SL(2,\C)$. The standard complex volume form on $\C^2$ induces via the isomorphism $\C^2\ra\Q$ and the inclusion $\Q\hookrightarrow\Im\Q\oplus\Q$ given above a complex 2-form $\theta$ on $\Im\Q\oplus\Q$. It is straightforward to verify that $\Phi_6^\ast\theta=-\gamma_3-i\gamma_2$. Hence, $\SL(2,\C)\subset\Stab(\gamma_2)\cap\Stab(\gamma_3)$. 
\end{proof}

We will need one more subgroup of $\GL(\mV)$. Put
\begin{equation}
\gR^+:=\bigg\{\varphi\in\GL(\mV)\bigg|\ \left(
\begin{matrix}
\lambda^{-2}&0&0\\
0&\lambda.1_2&0\\
0&0&\lambda^{-\frac{1}{2}}.1_4\\
 \end{matrix}
\right),\ \lambda>0
\bigg\}.
\end{equation}
It is easy to check that $\gR^+$ is a subgroup of $\gO_6$. From Lemma \ref{lemma inclusion type 6} follows that $\gR^+$ commutes with $\gH_6$. 
By the summary given above, $g.m. g^{-1}\in\gM$ whenever $m\in\gM,\ g\in\gH_6\times\gR^+$.  This shows that $\gO_6$ contains a subgroup $\gM\rtimes(\gH_6\times\gR^+)$.

\begin{lemma}\label{lemma semidirect product type 6}
$\gO_6=\gM\rtimes(\gH_6\times\gR^+)$.
\end{lemma}
\begin{proof}
Put $\gO_6':=\gM\rtimes(\gH_6\times\gR^+)$. Then $\mu|_{\gO_6'}$ induces epimorphism $\gO_6'\ra\CO(Q_2)$. By Lemma \ref{lemma split ses with H_6'}, we get a split short exact sequence $0\ra\gL'\ra\gO_6'\ra\CO(Q_2)\ra0$ where  $\gL'=\gM\rtimes\Phi_6^\ast\SL(2,\C)$. Then $\rho|_{\gL'}$ induces a split short exact sequence $0\ra\gM\ra\gL'\ra\Phi_6^\ast\SL(2,\C)\ra0$. By the proof Lemma \ref{lemma inclusion type 6}, $\Phi_6^\ast\SL(2,\C)\subset\Stab(\gamma_2)\cap\Stab(\gamma_3)$. It is straightforward to verify that also the other inclusion $"\supset"$ holds. The 2-form $\gamma_{i},\ i=2,3$ descends to the form $\sigma_{i-1}$ on $\mZ_4$. Thus $\Phi_6^\ast\SL(2,\C)=\Stab(\sigma_1)\cap\Stab(\sigma_2)$ and so $\Phi_6^\ast\SL(2,\C)$  is a splitting of the subgroup $\SL(2,\C)$ from (\ref{semidirect product type 6}). The claim now follows from the definition of $\gM$.
\end{proof}

\begin{thm}\label{thm max compact type 6}
 $\gK_6$ is a maximal compact subgroup of $\gO_6$.
\end{thm}
\begin{proof}
By Lemma \ref{lemma max compact in upper trian}, $\{1_7\}$ is a maximal compact subgroup of $\gM$. The same is obviously true also for $\gR^+$. From Lemma \ref{lemma split ses with H_6'} follows that $\gK_6$ is a maximal compact subgroup of $\gH_6$. The claim is then a consequence of Lemma \ref{lemma split ses with H_6'}.
\end{proof}

\subsection{Type 3}\label{section type 3}\label{section type 3}
A representative is
\begin{equation}
\omega_3=\alpha_{123}-\alpha_{167}+\alpha_{145}.
\end{equation}
Then $\mV_6:=\triangle_3^2=\triangle_3^3=[e_2,\ldots,e_7]$. This is an $\gO_3$-invariant subspace. We put $\mW_1:=\mV/\mV_6$. The insertion map $v\mapsto i_v\omega_3$ induces a monomorphism $\lambda:\mW_1\ra\Lambda^2\mV_6^\ast$. The image of $\lambda$ is an $\gO_3$-invariant conformally symplectic structure on $\mV_6$. The restriction map $\varphi\mapsto\varphi|_{\mV_6}$ induces a split short exact sequence $0\ra\gK\ra\gO_3\ra\CSp(3,\R)\ra0$. Moreover, $\gK= 
\bigg\{\varphi\in\GL(\mV)
\bigg|
\ \varphi_{\cB_{st}}=
\left(
\begin{matrix}
 1&\ast\\
 0&1_6
\end{matrix}
\right)
\bigg\}.$
\medskip

\begin{center}
 \textemdash\ x \textemdash
\end{center}

We put $f_1:=e_1,\ f_2:=e_2,\ f_3:=e_4,\ f_4:=e_6,\ f_5:=e_7,\ f_6:=e_5,\ f_7:=e_3$. Let $\cB_3^\ast=\{\beta_1,\dots,\beta_7\}$ be a dual basis to  $\cB_3:=\{f_1,\dots,f_7\}$. Then
\begin{equation}
 \omega_3=\beta_{125}+\beta_{136}+\beta_{147}.
\end{equation}

\bigskip

Let $\omega$ be the imaginary part of the standard Hermitian form on $\C^3$. Let  $\CSp(\omega):=\Stab([\omega])$ and $\Sp(\omega):=\Stab(\omega)$. The map $\mu:\CSp(\omega)\ra\R^\ast$ determined by $ \mu(\varphi).\omega=\varphi.\omega$ is a group homomorphism. 

Let $C:\C^3\ra\C^3$ be the standard conjugation. Then $C.\omega=-\omega$ and so $C\in\CSp(\omega),\ \mu(C)=-1$. Let $\wt\gU(3)$ be the group generated by $\gU(3)$ and $C$. It is easy to see that $\wt\gU(3)$ has two connected components, the connected component of the identity $\wt\gU(3)_0=\gU(3)$ and  $\wt\gU(3)_1:=\{\varphi\circ C:\ \varphi\in\gU(3)\}$.

\begin{lemma}\label{lemma split ses with CSp}
There is a split short exact sequence
\begin{equation}
 0\ra\Sp(\omega)\ra\CSp(\omega)\xrightarrow{\mu}\R^\ast\ra0.
\end{equation}
The group $\wt\gU(3)$ is a maximal compact subgroup of $\CSp(\omega)$.
\end{lemma}
\begin{proof}
Consider the subgroup $\gR^\ast$ of $\CSp(\omega)$ generated by $C$ and $\{\lambda.Id_{\C^3}:\ \lambda>0\}$. Then it is easy to see that $\mu$ restricts to an isomorphism $\gR^\ast\ra\R^\ast$. The inverse of this is a splitting of $\mu$.

$\{Id_{\C^3},C\}=\gO(\C^3)\cap\gR^\ast$ and $\gU(3)=\Sp(\omega)\cap\gO(3)$. As $\CSp(\omega)$ and $\gR^\ast$ are closed in $\GL(\C^3)$, the second claim is a corollary of Lemma \ref{lemma max compact in product}.
\end{proof}

For $A\in\CSp(\omega)$  we define $\phi_A\in\GL(\R\oplus\C^3),\ \phi_A(x,v)=(\mu(A).x,Av)$.  Let  
\begin{equation}
\Phi_3:\mV\ra\R\oplus\C^3,\ \
\Phi_3\big(\sum_{i=1}^7 x_if_i\big)=(x_1,(x_2+ix_5,x_2+ix_6,x_4+ix_6)).
\end{equation}
We put $\gH_3:=\{\Phi_3^\ast(\phi_A):\ A\in\CSp(\omega)\}$ and $\gK_3:=\{\Phi_3^\ast(\phi_A):\ A\in\wt\gU(3)\}$.

\begin{lemma}
 $\gH_3\subset\gO_3$.
\end{lemma}
\begin{proof}
Let $A\in\CSp(\omega)$ and  $\alpha\in\R^\ast$ be the linear form corresponding to $Id_\R$. Then $\alpha\wedge\omega\in\Lambda^3(\R\oplus\C^3)^\ast$ and $$\phi_A.(\alpha\wedge\omega)=\mu(A)^{-1}.\alpha\wedge A.\omega=\mu(A)^{-1}.\alpha\wedge\mu(A).\omega=\alpha\wedge\omega.$$ Since $\omega_3=\Phi_3^\ast(\alpha\wedge\omega)$, the claim follows.
\end{proof}

It is clear that $\gK$ and $\gH_3$ generate a subgroup $\gK\rtimes\gH_3$ of $\gO_3$.

\begin{lemma}
 $\gO_3=\gK\rtimes\gH_3$.
\end{lemma}
\begin{proof}
The restriction map $\varphi\in\gO_3\mapsto\varphi|_{\mV_3}$ induces isomorphism $\gH_3\ra\CSp(3,\R)$. This readily proves the claim.
\end{proof}

\begin{thm}\label{thm max compact type 3}
The group $\gK_3$  is a maximal compact subgroup of $\gO_3$.
\end{thm}
\begin{proof}
From  Lemma \ref{lemma max compact in upper trian} follows that  $\{1_7\}$ is a maximal compact subgroup of $\gK$. The claim is a consequence of  Lemmata \ref{lemma split ses with CSp} and \ref{lemma max compact in product I}.
\end{proof}

It is clear  that $\gK_3\subset\SO(7)$. We will now show $\gK_3$ can be viewed as a subgroup of $\Spin^c(7)$,  i.e. there is an embedding $\gK_3\ra\Spin^c(7)$ such that the composition $\gK_3\ra\Spin^c(7)\ra\SO(7)$ induces the identity on $\gK_3$. Recall that $\Spin^c(n)$ is the quotient of $\Spin(n)\times\gU(1)$  by the subgroup $\{\pm(1,1)\}$ where $\{\pm1\}=\ker\rho_n$ and  $\rho_n:\Spin(n)\ra\SO(n)$ is the usual 2:1 covering. We denote the class of $(a,e^{it})$ in the quotient  by $\langle a,e^{it}\rangle$. The canonical homomorphism $\rho_n^c:\Spin^c(n)\ra\SO(n)$ is  $\langle a,e^{it}\rangle\mapsto\rho_n(a)$.

\begin{lemma}\label{lemma subgroup of Spin^c(7)}
$\gK_3$ is a subgroup of $\Spin^c(7)$.
\end{lemma}
\begin{proof}
$\wt\gU(3)$ has two connected components and so the same is true for $\gK_3$.
The connected component $\gK_3^o$  of the neutral element of $\gK_3$ is isomorphic to $\gU(3)$. It is well known (see  \cite[Section 3.4]{Mo}) that $\gU(3)$ is a subgroup of $\Spin^c(6)$. Using  the standard inclusion $\Spin^c(6)\hookrightarrow\Spin^c(7)$, we see that there is  a subgroup $\wt\gK_3^o$  of $\Spin^c(7)$ such that  $\rho_7^c|_{\wt\gK_3^o}$ induces isomorphism   of Lie groups $\wt\gK_3^o\ra\gK_3^o$.  Hence, it remains to show that there is a subgroup  $\tilde\gK_3$ of $\Spin^c(7)$ which contains $\tilde\gK_3^o$ such that $\rho_7^c$ restricts to an isomorphism of Lie groups $\tilde\gK_3\ra\gK_3$.

Notice that $\psi_{C}:=\Phi_3^\ast(\phi_A)\in\SO(7)$ is determined by $f_i\mapsto- f_i,,\ i=1,2,3,4$ and $\ f_{j}\mapsto f_j,\ j=5,6,7$. Let $B_{\mV}$ be the standard inner product on $\mV$. Let us view $\Spin(7)$ as the subgroup of the Clifford algebra $(\mV,B_{\mV})$ that is generated by even number of unit vectors from $\mV$.  Put  $\alpha:= f_1.f_3.f_5.f_7\in\Spin(7)$. Then $\psi_C=\rho_7(\alpha)$ and $\alpha^2=1$. It is easy to verify that  we can take as $\tilde\gK_3$ the subgroup of $\Spin^c(7)$ that is generated by $\tilde\gK_3^o$ and $\langle\alpha,1\rangle$.
\end{proof}

\subsection{Type 4}\label{section type 4}
A representative of the orbit is
\begin{equation}
\omega_4=\alpha_{123}-\alpha_{167}+\alpha_{145}+\alpha_{246}.
\end{equation}
With respect to the basis $\cB_3$ from  Section \ref{section type 3} we have that 
\begin{equation}
\omega_4=\beta_{125}+\beta_{136}+\beta_{147}+\beta_{234}.
\end{equation}

Put $\mV_3:=[f_5,f_6,f_7],\ \mV_6:=[f_2,f_3,f_4,f_5,f_6,f_7]$. It is claimed in \cite[Section type 3]{BV} that $\triangle^2_4=\mV_3$ and $\triangle_4^3=\mV_6$ and so these are $\gO_4$-invariant subspaces. We put $\mW_1:=\mV/\mV_6, \mW_4:=\mV/\mV_3$ and $\mZ_3:=\mV_6/\mV_3$.

Each $\varphi\in\gO_4$ descends to  $\p\varphi\in\GL(\mW_1)\cong\R^\ast$ and $\hat\varphi\in\GL(\mZ_3)$. The homomorphism $\mu:\gO_4\ra\GL(\mW_1),\ \varphi\mapsto\p\varphi$  induces a split short exact sequence $0\ra\gO_4^+\ra\gO_4\xrightarrow{\mu}\R^\ast\ra0$. The 3-form $\omega_4$ induces a volume form on $\mZ_3$. It follows that the image of the homomorphism $\nu:\gO_4\ra\GL(Z_3),\ \varphi\mapsto\hat\varphi$ is contained in $\SL(\mZ_3)$. There is  a split short exact sequence $0\ra\gK\ra\gO_4^+\ra\SL(Z_3)\ra0$. If $\varphi\in\gK$, then $\varphi_{\cB_3}=
\left(
\begin{matrix}
1&0&0\\
\ast&1_3&0\\
\ast&\ast&1_3\\
\end{matrix}
\right)$. Moreover, it can be shown that $\gK\cong(((\R^3\rtimes\R^6)\rtimes\R)\rtimes\R)\rtimes\R$.
Fixing an isomorphism $\mZ^3\ra\R^3$, we obtain an isomorphism between $\gO_4$ and a semi-direct product
\begin{equation}
 (((((\R^3\rtimes\R^6)\rtimes\R)\rtimes\R)\rtimes\R)\rtimes\SL(3,\R))\rtimes\R^\ast
\end{equation}
as in \cite[7. Proposition]{BV}.

\begin{center}
 \textemdash\ x \textemdash
\end{center}

Put
\begin{equation}
\gH_4:=\Bigg\{\varphi\in\GL(\mV):  \varphi_{\cB_3}=\left(
 \begin{matrix}
1&0&0\\
0&A&0\\
0&0&(A^T)^{-1}\\
 \end{matrix}
 \right),\ A\in\SL(3,\R) \Bigg\}.
\end{equation}
It is clear that  $\gK_4^o:=\gH_4\cap\SO(7)$  is a maximal compact subgroup of $\gH_4$.

\begin{lemma}
$\gH_4\subset\gO_4$. 
\end{lemma}
\begin{proof}
Notice that $\omega_4=\omega_3+\beta_{234}$. It is obvious that $\gH_4\subset\Stab(\beta_{234})\cap\Stab(\beta_{1})$ and  it is straightforward to verify that $\gH_4\subset\Stab(\beta_{25}+\beta_{26}+\beta_{37})$. 
\end{proof}

Put 
\begin{equation}\label{R subgroup of O4}
\gR^\ast=\Bigg\{\varphi\in\GL(\mV):\ \varphi_M=\left(
 \begin{matrix}
\lambda&0&0\\
0&1_3&0\\
0&0&\lambda^{-1}.1_3\\
 \end{matrix}
 \right),\ \lambda\in\R^\ast\Bigg\}.
\end{equation}
It clear that $\gR^\ast\subset\gO_4$ and that $\gR^\ast$ commutes with $\gH_4$. Given $k\in\gK,\ g\in\gH_4\times\gR^\ast$, then we easily check that $g.k. g^{-1}\in\gK$. This shows that $\gO_4$ contains a subgroup $\gK\rtimes(\gH_4\times\gR^\ast)$.

\begin{lemma}
  $\gO_4=\gK\rtimes(\gH_4\times\gR^\ast)$.
\end{lemma}
\begin{proof}
 This easily follows from the summary given above.
\end{proof}

$\gR^\ast\cap\gO(7)$ contains two elements that correspond to $\lambda=\pm1$. We denote the intersection by $\Z_2$.

\begin{thm}\label{thm max compact type 4}
The group $\gK_4:=\gK^o_2\times\Z_2$ is a  maximal compact subgroup of $\gO_4$.  The group  $\gK_4$ is a subgroup of $\G_2$.
\end{thm}
\begin{proof}
By  Lemma \ref{lemma max compact in upper trian}, $\{1_7\}$ is maximal compact subgroup of $\gK$. Hence, it is (see Lemma \ref{lemma max compact in product I}) enough to show that $\gK_4$ is a maximal compact subgroup of $\gH_4\times\gR^\ast$. We have that $\gK_4=(\gH_4\times\gR^\ast)\cap\gO(7),\ \gK_4^o=\gH_4\cap\gO(7)$ and $\Z_2=\gR^\ast\cap\gO(7)$. Since $\gH_4$ and $\gR^\ast$ are closed subgroups of  $\GL(\mV)$ and since $\gK_4^o$ and $\Z_2$ is a maximal compact subgroup of $\gH_2$ and $\gR^\ast$, respectively, the claim follows from Lemma \ref{lemma max compact in product}.
 
To complete the proof, let  $\Phi_4:\mV\ra\Im\Q\oplus\Q$ be the linear isomorphism that sends the basis $\cB_3$ to the basis $$\{(0,1),(i,0),(j,0),(k,0),(0,i),(0,j),(0,k)\}.$$ Then it is easy to see that $\gK_4=\Phi_4^\ast\{\phi_{p,\pm 1.p}:\ p\in\Sp(1)\}$. Thus, $\gK_4\subset\Phi_4^\ast\G_2$
\end{proof}

\subsection{Type 1}\label{section type 1}
A representative of the orbit is
\begin{eqnarray}
\omega_1=\alpha_{127}+\alpha_{134}+\alpha_{256}.
\end{eqnarray}
 We have $\triangle^2_1=\mV^a_3\cup\mV^b_3$ and $\triangle_1^3=\mV_6^a\cup\mV_6^b$ where
\begin{equation*}
 \mV^a_3:=[e_3,e_4,e_7],\ \mV^b_3:=[e_5,e_6,e_7], \ \mV_6^a=[e_1,e_3,\dots,e_7],\ \mV_6^b=[e_2,\dots,e_7].
\end{equation*}
A subspace  $\mV_1:=\mV^a_3\cap\mV^b_3$ is $\gO_1$-invariant and we put $\mZ_2^a:=\mV_3^a/\mV_1,\ \mZ^b_2:=\mV^b_3/\mV_1$.
Each element $\varphi\in\gO_1$ induces an automorphism $\p\varphi$ of $\mZ_2^a\oplus\mZ_2^b$ such that $\p\varphi(\mZ_2^a)=\mZ_2^a,\ \p\varphi(\mZ_2^b)=\mZ_2^b$ or $\p\varphi(\mZ_2^a)=\mZ_2^b,\ \p\varphi(\mZ_2^b)=\mZ_2^a$. A map 
\begin{displaymath}
 sg:\gO_1\ra\Z_2,\ sg(\varphi)=
\Big\{
\begin{matrix}
1&\  \mathrm{if}\ \ \p\varphi(\mZ_2^a)=\mZ_2^a,\ \p\varphi(\mZ_2^b)=\mZ_2^b\\
-1&\   \mathrm{if}\ \ \p\varphi(\mZ_2^a)=\mZ_2^b,\ \p\varphi(\mZ_2^b)=\mZ_2^a\\
\end{matrix}
\end{displaymath}
is a split group homomorphism. Put $\gO_1^+:=\ker(sg)$. If $\varphi\in\gO_1^+$, then it is natural to view $\tilde\varphi$ as an element of $\GL(\mZ_2^a)\times\GL(\mZ_2^b)$. A map $\varphi\mapsto\p\varphi$ induces  a split short exact sequence $0\ra\gK\ra\gO_1^+\ra\GL(\mZ_2^a)\times\GL(\mZ_2^b)\ra0$. Moreover, it can be shown that
$\varphi_\cB=
\left(
\begin{matrix}
1_2&0&0\\
\ast&1_4&0\\
\ast&\ast&1
\end{matrix}
\right)$ for every  $\varphi\in\gK$ and that $\gK\cong(H\oplus H)\rtimes\R^4$ where $H$ is a Heisenberg algebra in dimension 3. All together, $\gO_1$ is isomorphic to a semi-direct product
\begin{equation}
 (((H\oplus H)\rtimes\R^4)\rtimes\GL(\mZ_2^a)\times\GL(\mZ_2^b))\rtimes\Z_2
\end{equation}
as in \cite[4. Proposition]{BV}.

\begin{center}
 \textemdash\ x \textemdash
\end{center}

We can now continue with the proof.
Let $A,B\in\GL(2,\R)$. We define 2-linear automorphisms of $\R\oplus\R\oplus\R^2\oplus\R^2\oplus\R$, namely
\begin{align*}
&\phi^+_{A,B}(u_1,u_2,v_1,v_2,u_3):=\bigg(\frac{u_1}{\det A},\frac{u_2}{\det B}, A(v_1),B(v_2),\det A.\det B.u_3\bigg)\ \mathrm{and}\\
&\phi^{-}_{A,B}(u_1,u_2,v_1,v_2,u_3):=\bigg(\frac{u_2}{\det B},\frac{u_1}{\det A}, B(v_2),A(v_1),-\det A.\det B.u_3\bigg)
\end{align*}
where $u_1,u_2,u_3\in\R,\ v_1,v_2\in\R^2$. Let 
$$\Phi_1:\R^7\ra \R\oplus\R\oplus\R^2\oplus\R^2\oplus\R,\ \Phi_1(x_1,\dots,x_7)=(x_1,x_2,(x_3,x_4),(x_5,x_6),x_7)$$
and put $\varphi^\bullet_{A,B}:=\Phi_1^\ast(\phi^\bullet_{A,B})$ where $\bullet=\pm$. 
It is straightforward to verify that 
\begin{eqnarray}
\varphi^+_{A,B}\circ\varphi^+_{C,D}=\varphi^+_{A.C,B.D},\  \varphi^{-}_{A,B}\circ\varphi^{-}_{C,D}=\varphi^+_{B.C,A.D},\\
\varphi^+_{A,B}\circ\varphi^{-}_{C,D}=\varphi^{-}_{B.C,A.D},\  \varphi^{-}_{A,B}\circ\varphi^{+}_{C,D}=\varphi^+_{C.A,D.B}
\end{eqnarray}
and so $\gH_1:=\{\varphi^\pm_{A,B}|\ A,B\in\GL(2,\R)\}$ is a subgroup of $\GL(\mV)$. The following observation is straightforward.

\begin{lemma}\label{lemma inclusion type 1}
$\gH_1\subset\gO_1$ 
\end{lemma}

It is  easy to check that $h.k.h^{-1}\in\gK$ whenever $h\in\gH_1,\ k\in\gK$. Hence, $\gO_1$ contains a subgroup $\gK\rtimes\gH_1$.

\begin{lemma}\label{lemma semi-direct product type 1}
 $\gO_1=\gK\rtimes\gH_1$
\end{lemma}
\begin{proof}
It is easy to see that the inclusion $\gH_1\hookrightarrow\gO_1$ is a splitting of the subgroup $(\GL(Z_2^a)\times\GL(Z_2^b))\rtimes\Z_2$. The claim follows from the definition of $\gK$.
\end{proof}

\begin{thm}\label{thm max compact type 1}
A maximal compact subgroup of  $\gH_1$ is a maximal compact subgroup of $\gO_1$.
\end{thm}
\begin{proof}
From Lemma \ref{lemma max compact in upper trian} follows that  $\{1_7\}$ is a maximal compact subgroup of $\gK$. The claim is then a consequence of Lemma \ref{lemma max compact in product I}.
\end{proof}

\section{Multisymplectic 3-forms of a fixed algebraic type}\label{section global forms}
\subsection{Characteristic classes of spin and spin$^c$ vector bundles}\label{section spin char class}
We will denote by $\underline\R^i$ a trivial vector bundle with fiber $\R^i$. For  a real vector bundle $\xi$ we denote by $w_i(\xi), \ p_1(\xi)$ and  $e(\xi)$ the $i$-th Stiefel-Whitney class, the first Pontryagin class and the Euler class of $\xi$, respectively. If $\xi$ is a complex vector bundle, then $c_i(\xi)$ is the $i$-th Chern class of $\xi$. We will need (see also \cite[Section 2]{CCV}) two more characteristic classes, one is defined for a spin vector bundle and the other is defined for a spin$^c$ vector bundle. 

Suppose that $\xi$ is a spin vector bundle over a base $N$. Then there is a class $q(\xi)\in H^4(N,\Z)$ which is independent of the choice of a spin structure on $\xi$. We have: 
\begin{equation}
 q(\xi\oplus\xi')=q(\xi)+ q(\xi') \ \ \mathrm{and} \ \ 2q(\xi)=p_1(\xi)
\end{equation}
if $\xi$ and $\xi'$ are spin. If $\xi$ is a complex bundle, then $\xi$ admits a spin structure if, and only if $c_1(\xi)$ is divisible by two, say $2m=c_1(\xi)$ for some $m\in H^2(M,\Z)$. Then 
\begin{equation}
 q(\xi)=2m^2-c_2(\xi).
\end{equation}
If $\xi=TN$, then we put $q(N):=q(TN)$.

A vector bundle $\xi$ admits a spin$^c$ structure if, and only if $w_2(\xi)$ is a reduction of an integral class, say $\rho_2(\ell)=w_2(\xi)$ for some $\ell\in H^2(N,\Z)$. Then we can define $q(\xi;\ell):=q(\xi-\lambda)\in H^4(N,\Z)$ where $\lambda$ is a complex line bundle with $c_1(\xi)=\ell$. The spin$^c$ characteristic class satisfies: 
\begin{equation}
2q(\xi;\ell)=p_1(\xi)-\ell^2,\ \rho_2(q(\xi;\ell))=w_4(\xi),\ q(\xi;l+2m)=q(\xi,\ell)-2\ell m-2m^2 
\end{equation}
where $m\in H^2(N,\Z)$. If $\xi$ is a complex vector bundle with $c_1(\xi)=\ell$, then
\begin{equation}
q(\xi;\ell)=-c_2(\xi).
\end{equation}
If $\xi=TN$, then  we put $q(N;\ell):=q(TN;\ell)$.

The following Theorem can be found in  \cite{CCS}. 

\begin{thm}\label{fundamental thm}
Let $N$ be a closed, connected manifold of dimension 7. Consider two orientable 7-dimensional real vector bundles $\xi$ and $\xi'$ over $N$ with $w_2(\xi) = w_2(\xi')=\rho_2(\ell)$, where $\ell \in H^2(N;\mathbb Z)$. Then $\xi$ and $\xi'$ are isomorphic as vector bundles if, and only if $q(\xi;\ell) = q(\xi'; \ell)$
\end{thm}

\subsection{Multisymplectic 3-forms of algebraic type 5,6,7,8}
\begin{thm}\label{thm global forms 5,6,7,8}
Let $N$ be a closed connected 7-dimensional manifold. Then the following are equivalent:
\begin{enumerate}
\item $N$ is orientable and spin.
\item $N$ admits a multisymplectic 3-form of algebraic type 8.
\item $N$ admits a multisymplectic 3-form of algebraic type 5.
\item $N$ admits a multisymplectic 3-form of algebraic type 6.
\item $N$ admits a multisymplectic 3-form of algebraic type 7.
\end{enumerate}
\end{thm}
\begin{proof}
(1)"$\Leftrightarrow$"(2) is proved in \cite{G}. 

(2)"$\Leftrightarrow$"(3) is proved in \cite{Le}. 

(1)"$\Rightarrow$"(4) and (5). 
By a result from \cite{Th}, any closed, orientable, spin manifold of dimension 7 admits two everywhere linearly independent vector fields. This gives a reduction from $\Spin(7)$-structure to $\Spin(5)$-structure. Following  \cite[Proposition 2.2]{CCS}, an inclusion 
$$\Sp(1)=\Sp(1)\times 1\hookrightarrow\Sp(1)\times\Sp(1)=\Spin(4)\hookrightarrow\Spin(5)$$ 
induces isomorphism on homotopy groups $\pi_i$ for $i\le 5$ and an epimorphism on $\pi_6$. Hence, the $\Spin(5)$-structure reduces to a $\Spin(4)$-structure which shows that $N$ admits three everywhere linearly independent vector fields, i.e. $TN\cong\underline\R^3\oplus\eta$ where $\eta$ is spin. Furthermore, the reduction from $\Spin(4)$ to $\Sp(1)$ means that $\eta$ has a $\Sp(1)$-structure or equivalently, $\eta$ is  a 1-dimensional $\Q$-vector bundle. Notice that a composition
$$\Sp(1)\hookrightarrow\Spin(5)\hookrightarrow\Spin(7)\xrightarrow{\rho_7}\SO(7)\ \mathrm{is}\footnote{Up to conjugation.}\ \ q\in\Sp(1)\mapsto\phi_{1,q}\in\SO(7).$$
This is an embedding $\Sp(1)\hookrightarrow\SO(7)$. We know from Section \ref{section type 6} that  $\Sp(1)\subset\gK_6$ and from Section \ref{section type 7} that $\Sp(1)\subset\gK_7$. We see that the $\Sp(1)$-structure extends to a $\gK_6$ and $\gK_7$-structure.

(4)"$\Rightarrow$" (3) By Section \ref{section type 6}, $\gK_6\subset\Phi_6^\ast\pG_2$.

(5)"$\Rightarrow$" (2) By  Section \ref{section type 7}, $\gK_7\subset\Phi_7^\ast\G_2$.
\end{proof}

\subsection{Multisymplectic 3-forms of algebraic type 4}
\begin{thm}\label{thm gl form 4}
Let $N$ be a closed, connected 7-manifold. If $N$ is orientable, spin and that there is $u\in H^4(M;\mathbb Z)$ such that $q(M) = -4u$, then 
$N$ admits a multisymplectic 3-form of algebraic type 4.

If $N$ admits a multisymplectic 3-form of algebraic type 4, then $N$ is orientable and spin.
\end{thm}
\begin{proof}
Recall Section \ref{section type 4} that the maximal compact, connected subgroup $\gK_4^o$ of $\gO_4$ is isomorphic to $\SO(3)$ and that $\mV\cong\R\oplus\R^3\oplus\R^3$ as a $\gK_4^o$-module where $\R$ is a trivial module and $\R^3$ is the standard representation. Alternatively, $\mV\cong\R\oplus\C^3$ where $\C^3=\R^3\otimes\C$. We see that $N$ admits a $\gK_4^o$-structure if, and only if there is a 3-dimension orientable vector bundle $\alpha$ over $N$ such that $TN\cong\underline\R\oplus\alpha_\C$ where  $\alpha_\C$ is the complexification of $\alpha$.

Let us assume that $q(N)=-4u$  for some $u\in H^4(N,\Z)$. We will first show that there is a $\Q$-line bundle $\mu$ such that  $-q(\mu)=e(\mu)=u$. For this we follow  \cite[Proposition 2.5]{CCS}. The map $e: BSU(2)\ra K(\Z,4)$ is an isomorphism on $\pi_i,\ i\le 4$ and epimorphism on $\pi_5$. This implies that there is  a $\Q$-line bundle $\eta$ over a 5-skeleton $N^{(5)}$ of $N$ with $e(\eta)=u$. Since $\pi_5(BO(\infty))=\pi_6(BO(\infty))=0$, it follows that the stable bundle $\eta\oplus\underline\R^3$ extends to $N$. Since any stable vector bundle over $N$ is stably isomorphic to a vector bundle of rank $7$, there is a  vector bundle $\xi$ over $N$ of rank 7 with $-q(\xi)=u$ and $w_2(\xi)=w_2(\mu)=0$. By a result from \cite{CS}, the bundle $\xi$ admits a $\Spin(5)$-structure. By the same argument as in the proof of Theorem \ref{thm global forms 5,6,7,8}, this structure reduces to $\Sp(1)$-structure, i.e. $\xi\cong\mu\oplus\underline\R^3$ where $\mu$ has a $\Sp(1)$-structure. We have  $-q(\mu)=-q(\xi)=u$.  

Next we take the associated 3-dimensional bundle $\alpha=\rho_-(\mu)$, see \cite[ Proposition 2.1]{CCVa}. From   \cite[Lemma 2.4]{CCVa}  follows that: $$p_1(\alpha)=p_1(\mu)-2e(\mu)=2q(\mu)-2e(\mu)=-4u=q(M).$$  On the other hand, $p_1(\alpha)=-c_2(\alpha_\C)=q(\alpha_\C)$. By Theorem \ref{fundamental thm}, $TN\cong\underline\R\oplus\alpha_\C$ and the sufficient condition follows.

By Theorem \ref{thm max compact type 4}, $\gK_4\subset\G_2$ and so the necessary condition follows from Theorem \ref{thm global forms 5,6,7,8}. 
\end{proof}

Notice that by  \cite[Lemma 2.6]{CD}, $0=w_4(M)=\rho_2(q(M))$ and so for a closed spin manifold $N$ there is always $v\in H^4(M;\mathbb Z)$ such that $q(M) = 2v$.

\subsection{Multisymplectic 3-forms of algebraic type 3}

\begin{thm}\label{global form 3}
Let $N$ be a 7-dimensional manifold. Then $N$ admits a multisymplectic 3-form of algebraic type 3 if, and only if $N$ is orientable and  spin$^c$.
\end{thm}
\begin{proof}
$"\Leftarrow"$ This is proved in \cite[Theorem 5.7]{D}. $"\Rightarrow"$
By Lemma \ref{lemma subgroup of Spin^c(7)}, $\gK_3$ is a subgroup of $\Spin^c(7)$ and so any $\gK_3$-structure extends to a $\Spin^c(7)$-structure. 
\end{proof}

\subsection{Multisymplectic 3-forms of algebraic type 2}

\begin{thm}\label{thm gl form 2}
Let $N$ be a 7-dimensional closed and connected manifold. If $N$ is orientable, spin and there are   $e,f\in H^2(M,\Z)$ such that  $e^2+f^2+3e f=-q(M)$, then $N$  admits a multisymplectic 3-form of algebraic type 2. If $N$ is simply-connected, then this condition is also necessary.

On the other hand, suppose that $N$ is orientable and admits a multisymplectic 3-form of algebraic type 2. Then $N$ is spin. 
\end{thm}
\begin{proof} Let $\alpha$ and $\beta$ be a complex line bundle with $c_1(\alpha)=e$ and $c_1(\beta)=f$. Put $\xi:=\underline\R\oplus\alpha\oplus\beta\oplus\alpha\otimes\beta$. By Theorem \ref{thm max compact type 2}, $\xi$ has a $\gK_2^o$-structure. On the other hand,  $w_2(\xi)=0$ and $q(\xi)=-c_2(\xi)=-e^2-f^2-3ef$. From Theorem \ref{fundamental thm} follows that $TN\cong\xi$. It is obvious that the condition is necessary, if  $N$ is simply-connected.

Let us assume that $N$ is orientable and admits a $\gO_2^+$-structure.
We have seen in the proof of Lemma \ref{lemma semi-direct product type 2}, that $\gO_2^+=\gO_2\cap\SL(\mV)=\gL\rtimes\gH_2$. Moreover, the inclusion $\gH_2\hookrightarrow\gO_2^+$ is a homotopy equivalence and thus any $\gO_2^+$-structure reduces to  $\gH_2$-structure. As $\gH_2\subset\Phi_6^\ast\pG_2$, any$\gH_2$-structure extends to a $\pG_2$-structure.  Theorem \ref{thm global forms 5,6,7,8} implies that $N$ is spin.  
\end{proof}

\subsection{Multisymplectic 3-forms of algebraic type 1}

\begin{lemma}
Let $N$ be a  7-manifold without boundary. Suppose that there are oriented vector bundles $\alpha,\beta$ of rank 2  such that $TN\cong\underline\R^3\oplus\alpha\oplus\beta$. Then $N$ admits a multisympletic 3-form of algebraic type 1. If $N$ is also simply-connected, then the assumption is also necessary.
\end{lemma}
\begin{proof}
Let us first show that the condition is sufficient. If $TN\cong\underline\R^3\oplus\alpha\oplus\beta$, then dually $T^\ast N\cong\underline\R^{3\ast}\oplus\alpha^\ast\oplus\beta^\ast$. Let  $\theta_i,\ i=1,2,3$ be differential forms on $N$ which trivialize $\underline\R^{3\ast}$.
Next we choose everywhere non-zero sections $\mu_1$ and $\mu_2$ of $\Lambda^2\alpha^\ast$ and $\Lambda^2\beta^\ast$, respectively. Now it is easy to see that $\theta_1\wedge\theta_2\wedge\theta_3+\theta_1\wedge\mu_1+\theta_2\wedge\mu_2$ is a multisymplectic 3-form of algebraic type 1.

Let us now assume that $\pi_1(N)=1$.  Recall  Lemma  \ref{lemma semi-direct product type 1} that $\gO_1=\gK\rtimes\gH_1$.  As $\gK\cong\R^{10}$ we know that $\pi_i(\gO_1/\gH_1)=1,\ i\ge0$. From the classical obstruction theory (see for example \cite{Th}) follows that any $\gO_1$-structure reduces to an $\gH_1$-structure. If $N$ is simply-connected, then any $\gH_1$-structure reduces to an $\gH_1^o$-structure where $\gH_1^o$  is  the connected component  of the identity element of $\gH_1$. Now $\gH_1^o\cong\GL^+(2,\R)^+\times\GL^+(2,\R)$ where $\GL^+(2,\R)=\{A\in\GL(2,\R):\ \det(A)>0\}$.  As an $\gH_1^0$-module, $\mV\cong\R^3\oplus\R^2_1\oplus\R^2_2$ where $\R^3$ is a trivial representation and $\R^2_i$ is the standard representation of the $i$-th factor of $\gH_1^o$. From this the necessary condition easily follows.
\end{proof}

We see that a simply-connected manifold $N$ with a global multisymplectic 3-form of algebraic type 1 admits a spin$^c$-structure. Proposition \ref{fundamental thm} implies the following.

\begin{thm}\label{thm gl form 1}
Let $N$ be a connected, closed and spin$^c$ 7-manifold. Suppose that there are $e,f\in H^2(M,\Z)$ such that $\rho_2(e+f)=w_2(M)$ and $-e f=q(M;e+f)$. Then $N$ admits a multisympletic 3-form of algebraic type 1. If $N$ is  simply-connected, then the assumption is also necessary.
\end{thm}

\end{document}